\documentclass[reqno]{amsart}
\usepackage[T1]{fontenc}
\usepackage{epsfig,amssymb,amsmath,mathrsfs,color}
\usepackage{hyperref}

\newcommand{\cB}{\mathcal{B}}
\newcommand{\cones}{\mathscr{T}}
\newcommand{\cK}{\mathcal{K}}
\newcommand{\R}{\mathbb{R}}
\newcommand{\N}{\mathbb{N}}

\newcommand{\zsets}{\mathcal{Z}}
\newcommand{\union}{\cup}
\newcommand{\intersection}{\cap}

\newcommand\rplus{\R_+}
\newcommand\euclidean{\partial C\backslash\{0\}}
\newcommand\euclideanT{\partial T\backslash[0]_T}
\newcommand\closure{\operatorname{cl}}
\newcommand\restsalt[2]{\closure\{g|_{#1} \mid g\in \cK_{\tau(#1,#2)} \}}

\newcommand\M[3]{\text{M}_{#3}(#1/#2)}
\newcommand\topfn{M}

\newcommand\hil[3]{\operatorname{Hil}_{#3}(#1,#2)}
\newcommand\funk[3]{\operatorname{Funk}_{#3}(#1,#2)}
\newcommand\Funky[1]{{\operatorname{Funk}_{#1} }}
\newcommand\revvy{\operatorname{R\!Funk}}
\newcommand\rev[3]{\operatorname{R\!Funk}_{#3}(#1,#2)}
\newcommand\minf{f}
\newcommand\h{f}
\newcommand\hfunk{f}
\newcommand\fun{h}
\newcommand\antiboundary{\mathcal{B}^\text{rev}}
\newcommand\minspacehilbert{\cB^\textrm{Hil}}
\newcommand\hrev{r}
\newcommand\hhil{h}

\newcommand\expf{j}
\newcommand\dotprod[2]{\langle{#1}|{#2}\rangle}
\newcommand{\conv}{\operatorname{conv}}
\newcommand\dualdbps{{}^*\!\mathcal{A}}
\newcommand\dbps{\mathcal{A}}
\newcommand\pul{\operatorname{Ls}}
\newcommand\pll{\operatorname{Li}}

\newtheorem{prop}{Proposition}[section]
\newtheorem{proposition}[prop]{Proposition}

\newtheorem{lemma}[prop]{Lemma}

\newtheorem{theorem}[prop]{Theorem}
\theoremstyle{definition}
\newtheorem{definition}[prop]{Definition}
\theoremstyle{remark}
\newtheorem{example}[prop]{Example}

\newcommand\mumax{\mu}
\newcommand\mumaxu{\mu_g}

\begin{document}

\title{The horofunction boundary of the Hilbert geometry}
\date{\today}
\author{Cormac Walsh}
\address{INRIA, Domaine de Voluceau,
78153 Le Chesnay C\'edex, France}
\email{cormac.walsh@inria.fr}


\begin{abstract}
We investigate the horofunction boundary of the Hilbert geometry defined on
an arbitrary finite-dimensional bounded convex domain $D$.
We determine its set of Busemann points, which are those points that
are the limits of ``almost-geodesics''.
In addition, we show that any sequence of points converging to a point
in the horofunction boundary also converges in the usual sense to a point
in the Euclidean boundary of $D$.
We prove that all horofunctions are Busemann points if and only if
the set of extreme sets of the polar of $D$ is closed in the
Painlev\'e--Kuratowski topology.
\end{abstract}

\maketitle

\section{Introduction}

\newcommand\vecspace{V}
\newcommand\scalar{\lambda}
\newcommand\linmorph{g}
\newcommand\interior{\mbox{int\,}}
\newcommand\poynta{x}
\newcommand\poyntb{y}
\newcommand\hilbertdist{\text{Hil}}

There has recently been growing interest in a certain metric space boundary,
which we call here the horofunction boundary.
To define this boundary for a metric space $(X,d)$, one assigns
to each point $z\in X$ the function $\phi_z:X\to \R$,
\begin{equation*}
\phi_z(x) := d(x,z)-d(b,z),
\end{equation*}
where $b$ is some basepoint.
If $X$ is proper, then the map
$\phi:X\to C(X),\, z\mapsto \phi_z$ defines an embedding of $X$ into $C(X)$,
the space of continuous real-valued functions on $X$ endowed
with the topology of uniform convergence on compact sets.
The horofunction boundary is defined to be
$X(\infty):=\closure\{\phi_z\mid z\in X\}\backslash\{\phi_z\mid z\in X\}$,
and its elements are called horofunctions.

This construction appears to be due to Gromov~\cite{gromov:hyperbolicmanifolds}.
Rieffel has an alternative construction based on
$C^*$-algebras~\cite{rieffel_group}. He calls it the metric boundary.

Of particular interest are those horofunctions that are the limits of
almost-geodesics. An almost-geodesic
as defined by Rieffel~\cite{rieffel_group}
is a map $\gamma$ from an unbounded set $T\subset \R_+$ containing 0 to $X$,
such that for any $\epsilon>0$,
\begin{equation*}
|d(\gamma(t),\gamma(s))+d(\gamma(s),\gamma(0))-t| < \epsilon
\end{equation*}
for all $t\in T$ and $s\in T$ large enough with $t\ge s$.
Rieffel calls the limits of such paths Busemann points.

The horofunction boundary is an additive version of the
Martin boundary appearing in probabilistic potential theory.
In~\cite{AGW-m}, this analogy was developed using max-plus (tropical) algebra.
The set of Busemann points was seen to be an
analogue of the \emph{minimal} Martin boundary.

There are few examples of metric spaces where the horofunction
boundary or Busemann points are explicitly known.
The first cases to be investigated
were those of Hadamard manifolds~\cite{ballmann:manifolds}
and Hadamard spaces~\cite{ballmann:spaces}, where the horofunction
boundary turns out to be homeomorphic to the ray boundary and all
horofunctions are Busemann points.
The case of finite-dimensional normed spaces has also received attention.
Andreev~\cite{andreev} found a connection with flag-directed sequences
and Karlsson~\emph{et.~al.}~\cite{karl_metz_nosk_horoballs}
determine the horofunction boundary when the norm is polyhedral.
A paper by the present author~\cite{walsh:normed} determines the
set of Busemann points of general finite-dimensional normed spaces.
Other examples of metric spaces where the horofunction
boundary has been studied include the Cayley graphs of
finitely-generated abelian groups, studied by
Develin~\cite{develin:compactification},
and Finsler $p$--metrics on $\text{GL}(n,\mathbb{C})/\text{U}_n$,
where explicit expressions for the horofunctions were found by
Friedland and Freitas~\cite{friedland_freitas_pmetrics,
friedland_freitas_revisiting1}.
Webster and Winchester have some general results on when all horofunctions
are Busemann points~\cite{winweb_busemann},~\cite{winweb_metric}.

In this paper we investigate the horofunction boundary of the Hilbert geometry.
Let $x$ and $y$ be distinct points in a bounded open convex subset
$D$ of $\R^N$, with $N\ge 1$.
Define $w$ and $z$ to be the points in the Euclidean boundary of $D$
such that $w$, $x$, $y$, and $z$ are collinear and arranged in this
order along the line in which they lie.
The Hilbert distance between $x$ and $y$ is
defined to be the logarithm of the cross ratio of these four points:
\begin{equation*}
\hil{x}{y}{}:= \log \frac{|zx|\,|wy|}{|zy|\,|wx|}.
\end{equation*}
If $D$ is the open unit disk, then the Hilbert metric is exactly
the Klein model of the hyperbolic plane.

As pointed out by Busemann~\cite[p.105]{busemann:geometry},
the Hilbert geometry is related to hyperbolic geometry in much the same way
that normed space geometry is related to Euclidean geometry.
It is not surprising therefore that there will be similarities between
the results here and those obtained in~\cite{walsh:normed} for normed spaces.

To state our main results it will be convenient to recall an alternative
definition of the Hilbert metric.
Let $C$ be an open cone in $\R^{N+1}$, which we take to mean
a non-empty subset of that space that is open, convex, invariant under
multiplication by positive scalars, and does not contain the origin.
Associated to $C$ is a relation $\le_C$ on $\R^{N+1}$ defined so that
$\poynta\le_C\poyntb$ if and only if $\poyntb-\poynta\in\closure C$.
If $C$ does not contain any lines, then this relation is a partial order.
In general it may fail to be anti--symmetric.

For each $\poynta\in C$ and $\poyntb\in \R^{N+1}$, define
\begin{align}
\label{eqn:funkcone}
\funk{\poyntb}{\poynta}{C}:=
   \log\inf\{\scalar >0 \mid \poyntb\le_C\scalar\poynta\}.
\end{align}
Then \emph{Hilbert's projective metric}
on the cone is defined to be
\begin{align}
\label{eqn:hilsplit}
\hil{\poynta}{\poyntb}{C}
   := \funk{\poynta}{\poyntb}{C}+\funk{\poyntb}{\poynta}{C},
\qquad\text{for all $\poynta$ and $\poyntb$ in $C$}.
\end{align}
If $C$ contains no lines, then
$\hilbertdist_C$ is a metric on the space of rays of the cone.
For further details, see the monograph of Nussbaum~\cite{nussbaum:hilbert}.
Other references for Hilbert's metric on a cone include~\cite{delaharpe},
\cite{gunawardenawalsh}, and~\cite{nussbaumwalsh}.

One can recover Hilbert's original definition, in the case when $C$
contains no lines, by taking $D$ to be
a cross section of $C$, that is letting $D:=\{ x\in C:\psi(x)=1 \}$,
where $\psi:\R^{N+1}\to \R$ is some linear functional that is positive
with respect to the partial ordering associated to $C$.
On $D$, which is a bounded convex open set, the two definitions agree.

Restricted to $D\times D$, the function $\Funky{C}$ can be written as
\begin{equation}
\label{eqn:funkratio}
\funk{x}{y}{C}:= \log \frac{|zx|}{|zy|},
\qquad\text{for all $x$ and $y$ in $D$.}
\end{equation}
Indeed, this is its usual definition.
Under this restriction, $\Funky{C}$ satisfies the usual metric
space axioms, apart from that of symmetry.
On $C$, it satisfies the triangle inequality
but can take negative values.

The expression of the Hilbert metric as the symmetrisation
of the Funk metric will play a crucial role in what follows.
It will turn out that every Hilbert horofunction is the sum of a
horofunction in the Funk geometry and a horofunction in the reverse Funk
geometry, where the metric in the latter is given by
\begin{equation*}
\rev{x}{y}{C}:=\funk{y}{x}{C}.
\end{equation*}
This will allow us to simplify the problem by investigating separately
the horofunction boundaries of these two geometries and then combining the
results. Determining the boundary of the Funk geometry turns out to be very
similar to determining that of a normed space,
which was done in~\cite{walsh:normed}.

We introduce a slight modification of the usual definition of
tangent cone, one more suited to dealing with open cones.
For any cone $T$ and point $x$ in its Euclidean boundary $\partial T$,
define the open tangent cone to $T$ at $x$ by
\begin{equation*}
\tau(T,x):= \{ \lambda(y-x)\mid \text{$\lambda>0$ and $y\in T$} \}.
\end{equation*}
Next, define the map between sets of cones that corresponds to taking all
open tangent cones of all members of the set:
\begin{equation*}
\Gamma(\mathbb{T})
   := \{ \tau(T,x) \mid  \text{$T\in \mathbb{T}$ and $x\in \partial T$} \},
\qquad\text{for any set of cones $\mathbb{T}$.}
\end{equation*}
Now iterate this map on any cone $T$ to get
\begin{equation*}
\cones(T) := \bigcup_{k=1}^\infty \Gamma^k(\{T\}).
\end{equation*}
Note that $T\in\cones(T)$.
For any open cone $T$ and point $p$ in $T$,
denote by $\hfunk_{T,p}$ the function from $T$ to~$\R$ defined by
\begin{equation*}
\hfunk_{T,p}(x) := \funk{x}{p}{T} - \funk{b}{p}{T},
\end{equation*}
where $b\in T$ is a basepoint.
Also define, for each $p\in \partial T$, the function $\hrev_{T,p}$
from $T$ to~$\R$ by
\begin{equation*}
\hrev_{T,p}(x) := \rev{x}{p}{T} - \rev{b}{p}{T}.
\end{equation*}
Note that in both cases the dependence on $p$ is only through the ray
on which $p$ lies.
%Observe that $\tau(T,x)$ contains $T$, and so if we have a basepoint
%in $T$, then we may reuse it as a basepoint in $\tau(T,x)$.
%The same is true for each cone in $\cones(T)$.

We use the notation $f|_S$ to denote the restriction of a function $f$
to a subset $S$ of its domain.

Now we can state our main results.
\begin{theorem}
\label{theorem1}
The set of Busemann points of the Hilbert geometry on a finite-dimensional
open cone $C$ containing no lines is
\begin{equation*}
\{ \hrev_{C,z} + \hfunk_{T,p}|_C \mid 
    \text{ $z\in \partial C\backslash\{0\}$,
           $T\in \cones(\tau(C,z))$,
           and $p\in T$
          } \}.
\end{equation*}
\end{theorem}
Although the expression in~(\ref{eqn:hilsplit}) for the Hilbert metric is
symmetric in $\Funky{}$ and $\revvy$, the description of the horofunctions
and Busemann points is not. This is because of the very different nature
of the boundaries of these two geometries. Some examples of horofunctions are
given in section~\ref{sec:examples}.

We use the theorem above to characterise those Hilbert geometries for which all
horofunctions are Busemann points.
Recall that a convex subset $E$ of a convex set $D$ is said to be an
\emph{extreme set} if the endpoints of any line segment in $D$ are
contained in $E$ whenever any interior point of the line segment is.
\begin{theorem}
\label{theorem3}
Let $D$ be a bounded convex open subset of $\R^N$ containing the origin.
Every horofunction of the Hilbert geometry on $D$ is a Busemann point
if and only if the set of extreme sets of the polar of $D$
is closed in the Painlev\'e--Kuratowski topology.
%A necessary and sufficient condition for every horofunction on a
%bounded convex open subset of $\R^N$ containing the origin
%to be a Busemann point in the Hilbert geometry is that
%the set of extreme sets of its polar
%be closed in the Painlev\'e--Kuratowski topology.
\end{theorem}

It had previously been shown~\cite{karl_metz_nosk_horoballs}
that all horofunctions of the Hilbert geometry on a polytope
are Busemann points.

In section~\ref{sec:examples}, we give examples of domains satisfying
and not satisfying the condition of Theorem~\ref{theorem3}.

Our final theorem generalises one of~\cite{karl_foertsch},
where the result was proved in the case when the sequence is a geodesic.

\begin{theorem}
\label{theorem2}
Let $D$ be a bounded convex open subset of $\R^N$.
If a sequence in $D$ converges to a point in the horofunction boundary of the
Hilbert geometry, then the sequence converges in the usual sense to a
point in the Euclidean boundary~$\partial D$.
\end{theorem}

So there is a continuous surjection from the horofunction compactification
to the Euclidean (usual) compactification that maps every point in $D$ to
itself. This is reminiscent of the continuous surjection from the horoboundary
of a $\delta$-hyperbolic space to its Gromov
boundary~\cite{coornaert_papadopoulos_horofunctions,winweb_hyperbolic,
storm_barycenter}.

\section{Boundary of the reverse-Funk geometry}

We start off by investigating the horofunctions and Busemann points of the
reverse Funk geometry.

A little care is needed here because $\revvy_T$ is not a metric and
so does not give rise to a topology. For this reason, we use the
topology of pointwise convergence instead of that of uniform convergence
on compact sets in the definition of the horofunction boundary of $\revvy_T$.

We will also find it convenient to use in the remainder of the paper
a slightly different definition of almost--geodesic, one adapted
from~\cite{AGW-m}.
Given on some space $X$ a function $d:X\times X\to \R$ satisfying the
triangle inequality, we say that a path $(x_l)_{l\in\N}$ is an
\emph{almost--geodesic} if, for some $\epsilon>0$,
\begin{align}
\label{eqn:def_almost}
d(x_0,x_1) + \dots + d(x_{l-1},x_l) \le d(x_0,x_l) + \epsilon,
\qquad
\text{for all $l\ge1$}.
\end{align}
We refer to $\epsilon$ as the parameter of the almost-geodesic.

We now define a Busemann point to be a pointwise limit of
$d(\cdot,x_l)-d(b,x_l)$ along an almost--geodesic that can not be written
as $d(\cdot,p)-d(b,p)$ for any $p\in X$.
It was shown in~\cite{AGW-m} that if $d(\cdot,\cdot)$ is a metric,
then this definition of Busemann point coincides with the one of Rieffel
discussed earlier.

We work in a finite-dimensional real vector space $V$
and denote its dual by $V^*$.
Recall that the dual cone $T^*$ of a cone $T$ is the set
\begin{align*}
\{z \in V^* \mid \text{$\dotprod{z}{x} \ge 0$ for all $x\in T$}\}.
\end{align*}
For any open cone $T$, we use the notation $[0]_T$ to denote the set
$\{x\in V \mid \text{$x\le_T 0$ and $0\le_T x$} \}$.
If $T$ contains no lines, then $[0]_T=\{0\}$.
For any open cone $T$, define
\begin{align}
\M{y}{x}{T}:=
   \inf\{\scalar >0 \mid y\le_T \scalar x\},
\qquad\text{for each $y\in V$ and $x\in T$.}
\end{align}
We extend the definition of $\revvy$ slightly by taking
$\rev{x}{y}{T}:=\log\M{y}{x}{T}$ for all $x\in T$ and $y\in V$.

We say that a sequence of points $(x_n)_{n\in\N}$ in $T$ converges in the
\emph{reverse Funk sense} to a function $f:T\to\R$ if
$\rev{\cdot}{x_n}{T} - \rev{b}{x_n}{T}$ converges pointwise to $f$ on $T$.
When we say that $(x_n)_{n\in\N}$ converges in the usual sense,
we mean with respect to the usual topology on $V$.

\begin{lemma}
\label{lem:dualformula}
Let $T\subset V$ be an open cone. Then
\begin{align*}
M_T(y/x)= \sup_{z\in T^*}\frac{\dotprod{z}{y}}{\dotprod{z}{x}},
\qquad
\text{for all $y\in V$ and $x\in T$.}
\end{align*}
\end{lemma}
\begin{proof}
We have that $\lambda x-y\in \closure T$ if and only if
$\lambda\dotprod{z}{x}-\dotprod{z}{y}\ge 0$ for all $z\in T^*$.
The conclusion follows.
\end{proof}

\begin{lemma}
\label{lem:continuity}
Let $T\subset V$ be an open cone.
The function $M_T(\cdot/\cdot)$ is jointly continuous in its two entries.
\end{lemma}
\begin{proof}
Since $T$ has non-empty interior, $T^*$ contains no lines, and so we can find
a compact cross section $D$ of $T^*$.
The expression $\dotprod{z}{y}/\dotprod{z}{x}$ remains unchanged when
$z$ is multiplied by a scalar and so, by Lemma~\ref{lem:dualformula},
we have
\begin{align*}
M_T(y/x)= \sup_{z\in D}\frac{\dotprod{z}{y}}{\dotprod{z}{x}},
\qquad
\text{for all $y\in V$ and $x\in T$.}
\end{align*}
As a supremum of a set of continuous functions, $M_T$ is
lower--semicontinuous.

Now let $(x_n)_{n\in\N}$ be a sequence in $T$ converging to a point $x\in T$
and let $(y_n)_{n\in\N}$ be a sequence in $V$ converging to $y\in V$.
Since the supremum is over a compact set, there is, for each $n\in\N$,
a point $z_n\in D$ such that
$M_T(y_n/x_n)=\dotprod{z_n}{y_n}/\dotprod{z_n}{x_n}$.
By taking a subsequence if necessary, we may assume that $M_T(y_n/x_n)$
converges to its limit supremum and furthermore, since $D$ is compact,
that the sequence $(z_n)_{n\in\N}$ converges to some $z\in D$.
So, the limit supremum of $M_T(y_n/x_n)$ is the limit of
$\dotprod{z_n}{y_n}/\dotprod{z_n}{x_n}$,
which is $\dotprod{z}{y}/\dotprod{z}{x}$.
However this is obviously no greater than $M_T(y/x)$.
We have thus proved that $M_T$ is upper--semicontinuous.
\end{proof}

\begin{lemma}
\label{lem:revlemmaB}
Let $T\subset V$ be an open cone and let $p\in \euclideanT$.
Let $z$ be in $T$ and define $y_\lambda := (1-\lambda)p+\lambda z$
for all $\lambda\in(0,1)$.
Then
\begin{equation*}
\lim_{\lambda\to 0}\hrev_{T,p}(y_\lambda) = -\rev{b}{p}{T}.
\end{equation*}
\end{lemma}
\begin{proof}
Let $B:=\sup\{\beta\in\R \mid \beta p \le_T z-p\}$.
Observe that $B<\infty$ since otherwise $p \le_T z/(\beta+1)$ for
arbitrarily large $\beta$, and this would imply that $p \le_T 0$,
contradicting our assumption on $p$.
Also, $0\le_T z$, and so $B\ge -1$. So $B$ is finite.
We have
\begin{align*}
\topfn_T(p/y_\lambda)
   &= \inf\{\alpha>0 \mid p \le_T \alpha((1-\lambda)p+\lambda z) \} \\
   &= \inf\{\alpha>0 \mid
                \Big(\frac{1}{\alpha} -1\Big)p \le_T \lambda(z-p) \} \\
   &= (\sup\{\beta>-1\mid \beta p \le_T \lambda(z-p) \} +1)^{-1} \\
   &= (\lambda\sup\{\gamma>-1 / \lambda \mid \gamma p \le_T z-p \} +1)^{-1}.
\end{align*}
As $\lambda$ tends to zero, the supremum in the last line converges
to $B$, and so $\topfn_T(p/y_\lambda)$ converges to $1$.
Therefore, $\rev{y_\lambda}{p}{T}$ converges to zero. The conclusion follows.
\end{proof}

\begin{lemma}
\label{lem:revconv}
Let $T\subset V$ be an open cone. If a sequence in $T$ converges in the
usual sense to a point $x$ in $\euclideanT$, then it converges in the
reverse Funk sense to $\hrev_{T,x}$.
\end{lemma}
\begin{proof}
Let $(x_n)_{n\in\N}$ be a sequence in $T$ converging to $x$.
Since $M_T$ is continuous by Lemma~\ref{lem:continuity}, $\rev{y}{\cdot}{T}$
is continuous for each $y\in T$. So $\rev{y}{x_n}{T}-\rev{b}{x_n}{T}$
converges to $\rev{y}{x}{T}-\rev{b}{x}{T} = \hrev_{T,x}(y)$
for each $y\in T$.
\end{proof}

\begin{proposition}
\label{pro:reversehorofunctions}
Let $C\subset V$ be an open cone containing no lines.
The set of horofunctions in the reverse Funk geometry on $C$ is
$\antiboundary:=\{\hrev_{C,x} \mid x\in\euclidean \}$.
All these  horofunctions are Busemann points of the reverse Funk geometry.
A sequence in a cross section of $C$ converges in the reverse Funk sense to
$\hrev_{C,x}\in \antiboundary$ if and only if it converges in the
usual sense to a positive multiple of $x$.
\end{proposition}
\begin{proof}
Let $(x_n)_{n\in\N}$ be a sequence in $C$ converging in the reverse Funk sense
to a horofunction $f$. We may assume that $(x_n)_{n\in\N}$ is
contained in some cross section $D$ of the cone having compact closure
and containing the base point $b$.
So $(x_n)_{n\in\N}$ has a limit point $x\in \closure D$ in the usual topology.
The point $x$ cannot be in $D$ for otherwise $f$ would not be a horofunction.
So $x\in\partial D$ and, by Lemma~\ref{lem:revconv}, $f=\hrev_{C,x}$.

It is well known~\cite[p.59]{shen:spray} that line segments are geodesics
in the Funk geometry, and hence also in the reverse Funk geometry.
Since every point in $\partial D$ is the endpoint of a line
segment contained in $D$, every horofunction in $\antiboundary$ is the limit
of a geodesic and therefore a Busemann point.

Now, let $x$ and $y$ be distinct points of $\partial D$.
For each $\lambda\in(0,1)$, let $x_\lambda:= (1-\lambda)x + \lambda b$
and $y_\lambda:= (1-\lambda)y + \lambda b$.
We calculate that
\begin{align*}
\lim_{\lambda\to 0} \hrev_{C,y}(x_\lambda)
     = \log \frac{|wy|}{|wx|} - \rev{b}{y}{C},
\end{align*}
where $w$ is the point in the intersection of the line $xy$ with
$\partial D$ the farthest from $y$ on the same side of $y$ as $x$.
Similarly,
\begin{equation*}
\lim_{\lambda\to 0} \hrev_{C,x}(y_\lambda)
    = \log\frac{|zx|}{|zy|} - \rev{b}{x}{C},
\end{equation*}
where $z$ is the point in the intersection of the line $xy$ with
$\partial D$ the farthest from $x$ on the same side of $x$ as $y$.
Also, by Lemma~\ref{lem:revlemmaB},
$\lim_{\lambda\to 0} \hrev_{C,x}(x_\lambda)=-\rev{b}{x}{C}$, and a
similar formula holds with $y$ instead of $x$.
Since both $\rev{b}{x}{C}$ and $\rev{b}{y}{C}$ are finite,
we deduce that
\begin{equation*}
\lim_{\lambda\to 0} (\hrev_{C,y}(x_\lambda) - \hrev_{C,x}(x_\lambda)
   + \hrev_{C,x}(y_\lambda) - \hrev_{C,y}(y_\lambda) )
       = \log\frac{|wy||zx|}{|wx||zy|}.
\end{equation*}
Since $x$ and $y$ are different, the right hand side is strictly positive,
which implies that $\hrev_{C,x}$ and $\hrev_{C,y}$ are different
since otherwise the left hand side would be zero.

Let $(z_n)_{n\in\N}$ be any sequence in $D$ and let $X$ be the set of
its limit points in the usual topology on $V$.
The set of its limit points in the
reverse Funk geometry will then be $\{\hrev_{C,x}\mid x\in X\}$
by Lemma~\ref{lem:revconv}.
From what we have
just proved, this set will contain a single point if and only if $X$
does. Thus, $(z_n)_{n\in\N}$ converges in the reverse Funk sense 
if and only if it converges in the usual sense.
\end{proof}

\section{Boundary of the Funk geometry}

To determine the Busemann points of the Funk geometry, it will be important
to investigate the relationship between the Funk geometry on a cone and that
on its tangent cones.

For the same reason as in the case of the reverse Funk geometry,
we define the horofunctions of the Funk geometry using pointwise convergence
and define almost--geodesics as in~(\ref{eqn:def_almost}).
We say that a sequence $(x_n)_{n\in\N}$ in an open cone $T$ converges to $f$
in the \emph{Funk sense} if $\funk{\cdot}{x_n}{T} - \funk{b}{x_n}{T}$
converges pointwise to $f$.
A Busemann point is the limit of an almost--geodesic in the Funk sense
not of the form $\funk{\cdot}{p}{T} - \funk{b}{p}{T}$ for some $p\in T$.

For each open cone $T\subset V$ and $x\in \partial T\backslash[0]_T$,
let $A_T(x)$ denote the set of Funk-geometry
horofunctions that may be attained as a limit in the Funk sense of a
sequence converging to $x$ in the usual sense.

\begin{lemma}
\label{lem:lemmaC}
Let $T$ be an open cone in $V$ and let $z$ be in $[0]_T$.
Then
\begin{align*}
\funk{(1-\alpha)z+\alpha x}{y}{T} &= \log \alpha + \funk{x}{y}{T}
\qquad\text{and} \\
\funk{x}{(1-\alpha)z+\alpha y}{T} &= -\log \alpha + \funk{x}{y}{T},
\end{align*}
for all $x,y\in T$ and $\alpha>0$.
\end{lemma}
\begin{proof}
Observe that $(1-\alpha)z+\alpha x\le_{T} \lambda y$ is equivalent to
$\alpha x\le_{T} \lambda y$, which in turn is equivalent to
$x\le_{T}\lambda y/\alpha$.
The first formula now follows on applying the definition of the Funk metric
in~(\ref{eqn:funkcone}). The proof of the second formula is similar.
\end{proof}

\begin{lemma}
\label{lem:closed}
Let $T$ be an open cone in $V$. Then,
$A_T(x)$ is closed for each $x\in \partial T\backslash[0]_T$.
\end{lemma}
\begin{proof}
Let $(g_n)_{n\in\N}$ be a sequence in $A_T(x)$ converging pointwise
to a Funk geometry horofunction $g$.
For each $n\in\N$, let $(x^{n}_i)_{i\in\N}$ be a sequence of points in $T$
converging in the usual sense to $x$ and in the Funk sense to $g_n$.
Also, let $(U_i)_{i\in\N}$ be some decreasing basis of open neighbourhoods
of $g$ in the horofunction compactification of the Funk geometry,
and $(V_i)_{i\in\N}$ be some decreasing basis of open
neighbourhoods of $x$ in the usual topology.
For each $i\in\N$, choose $n_i\in\N$ large enough that $g_{n_i}\in U_i$.
Having done this, we may choose $m_i$ large enough that
$x^{n_i}_{m_i}$ is in $V_i$ and
$\funk{\cdot}{x^{n_i}_{m_i}}{T}-\funk{b}{{x^{n_i}_{m_i}}}{T}$ is in $U_i$. So
$(x^{n_i}_{m_i})_{i\in\N}$ converges to $x$ in the usual sense and to $g$
in the Funk sense. Hence, $g\in A_T(x)$.
\end{proof}

For any open cone $T$ in $V$, define
\begin{align*}
\cK_T &:= \{ \minf_{T,p} \mid p\in T \}
\qquad\text{and}\\
\cB_T  &:= \{ f\in (\closure \cK_T)\backslash \cK_T
                \mid \text{$f$ is a Busemann point of $\Funky{T}$} \}.
\end{align*}
We use the notation $f|_X$ to denote the restriction of a function $f$ to
a set $X$.

\begin{figure}
\input{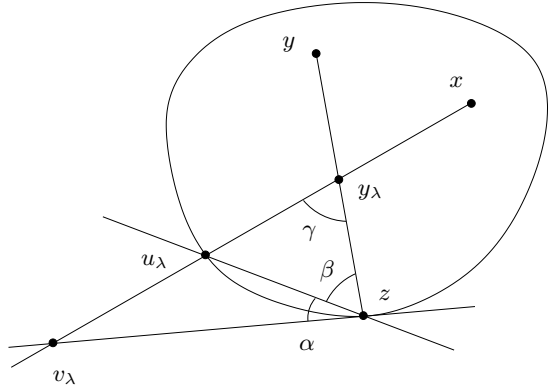}
\caption{Diagram for the proof of Lemma~\ref{lem:lemmaB}.
A two-dimensional cross section of the cone is shown.}
\label{fig:lemmaB}
\end{figure}

\begin{lemma}
\label{lem:lemmaB}
Let $x$ and $y$ be points of an open cone $T$ and let
$z$ be in the Euclidean boundary of $T$.
Define $y_\lambda:=(1-\lambda)z+\lambda y$ for all $\lambda\in(0,1)$.
Then
\begin{equation}
\label{eqn:tozero}
\funk{x}{y_\lambda}{T} - \funk{x}{y_\lambda}{\tau(T,z)} \to 0
\qquad
\text{as $\lambda\searrow 0$.}
\end{equation}
Moreover, $y_\lambda$ converges to $\minf_{\tau(T,z),y}|_T$
in the $\Funky{T}$ sense as $\lambda\searrow 0$.
Finally,
\begin{equation*}
\restsalt{T}{z} \subset \closure\cK_T \intersection A_T(z).
\end{equation*}
\end{lemma}
\begin{proof}
For $\lambda$ small enough, $y_\lambda$ is not greater than
$x$ in either the ordering associated to the cone $T$ or that associated to
$\tau(T,z)$. Therefore, for $\lambda$ small enough,
\begin{equation*}
\funk{x}{y_\lambda}{T} = \log\frac{|x u_\lambda|}{|y_\lambda u_\lambda|}
\qquad
\text{and}
\qquad
\funk{x}{y_\lambda}{\tau(T,z)}
   = \log\frac{|x v_\lambda|}{|y_\lambda v_\lambda|},
\end{equation*}
where $u_\lambda$ is the intersection with the boundary of $T$
of the half-ray starting at $x$ and passing through $y_\lambda$,
and $v_\lambda$ is the intersection of the same
half-ray with the boundary of $\tau(T,z)$.
Define the angles
\begin{equation*}
\alpha:=\angle u_\lambda z v_\lambda,
\qquad
\beta :=\angle y_\lambda z u_\lambda,
\qquad
\text{and}
\qquad
\gamma :=\angle v_\lambda y_\lambda z.
\end{equation*}
as in Figure~\ref{fig:lemmaB}.
Observe that $\alpha\to 0$ and $\gamma \to \angle yzx$ as $\lambda\to 0$,
whereas $\alpha+\beta$ remains constant.
So, applying the sine rule twice, we see that
\begin{equation*}
\frac{|y_\lambda u_\lambda|}{|y_\lambda v_\lambda|}
   = \frac{\sin \beta \sin(\pi-\alpha-\beta-\gamma)}
       {\sin(\alpha+\beta)\sin(\pi-\beta-\gamma)}
\longrightarrow 1
\qquad
\text{as $\lambda\searrow 0$.}
\end{equation*}
From this and the fact that both $|xu_\lambda|$ and $|xv_\lambda|$
converge to $|xz|$, the first statement of the lemma follows.

To prove the second, we apply~(\ref{eqn:tozero}) twice, first with $x$ an
arbitrary point in $T$ and then with $x=b$. We then use Lemma~\ref{lem:lemmaC},
which is applicable since $z$ is in $[0]_{\tau(T,z)}$.
The result is that $f_{T, y_\lambda}(x)$ converges to $f_{\tau(T,z),y}(x)$
for all $x\in T$.
This proves the second statement.

We have demonstrated that $\minf_{\tau(T,z),y}|_T$ is in $\closure \cK_T$
for all $y\in T$. But it is not difficult to show that,
for any $w\in\tau(T,z)$, there exists $y\in T$ such that
$w=(1-\lambda)z+\lambda y$ for some $\lambda>0$.
By Lemma~\ref{lem:lemmaC}, $f_{\tau(T,z),w}=f_{\tau(T,z),y}$.
Therefore the set $\{h|_{T} \mid h\in \cK_{\tau(T,z)} \}$
is a subset of $\closure\cK_T$.
By the second part of the present lemma, this set is also a subset of $A_T(z)$.
The third statement follows from this on taking closures since $A_T(z)$
is closed by Lemma~\ref{lem:closed}.
\end{proof}

In the next lemma, we will need the following notions.
Given on some space $X$ a function $d:X\times X\to \R$ satisfying the
triangle inequality and a function $g:X\to\R$,
we say that a path $(x_l)_{l\in\N}$ is an
\emph{almost--optimal path} with respect to $g$ if for some $\epsilon>0$,
\begin{equation*}
g(x_0) \ge -\epsilon + d(x_0,x_1) + \dots + d(x_{l-1},x_l) + g(x_l),
\qquad
\text{for all $l\ge1$}.
\end{equation*}

A \emph{min--plus measure} is a lower semicontinuous function from
some set to $\R\union\{+\infty\}$.
Let $g:T\to\R$ be a function.
The following map $\mumaxu: \closure \cK_T\to \R\union\{+\infty\}$
is a min--plus measure:
\begin{align*}
\mumaxu(w) := \inf \liminf_{x_n \to w}( \funk{b}{x_n}{T} + g(x_n)),
\qquad\text{for all $w\in \closure\cK_T$} ,
\end{align*}
where the infimum is taken over all sequences $(x_n)_{n\in\N}$
converging to $w$ in the Funk sense.

\begin{lemma}
\label{lem:representation}
Let $T$ be an open cone.
If $h \in A_T(x)$ with $x\in \partial T\backslash[0]_T$, then
\begin{equation}
\label{eqn:representation}
h(y) = \inf_{w\in\cB_T} \Big( w(y) + \nu(w) \Big),
\qquad\text{for all $y\in T$},
\end{equation}
for some min-plus measure $\nu$ on $\cB_T$ taking the value $+\infty$
outside $\restsalt{T}{x}$.
\end{lemma}
\begin{proof}
Let $y\in T$ and consider the straight line segment from $y$ to $x$.
By Lemma~\ref{lem:lemmaC}, the $\Funky{T}$ distance from $y$ to a point
increases in an unbounded fashion as the point moves along this line segment
from $y$ to $x$. So we may parameterise this line segment in such a way
as to get a map $\gamma:\rplus\to T$ such that $\gamma(0):= y$ and
$\funk{y}{\gamma(t)}{T}=t$ for all $t\in\rplus$.
By Lemma~\ref{lem:lemmaB}, $\gamma(t)$
converges to $f_{R,y}|_T$ in the $\Funky{T}$ sense as $t\to \infty$,
where $R:=\tau(T,x)$.

Let $(x_n)_{n\in\N}$ be a sequence in $T$ converging to $h$ in the
$\Funky{T}$ sense and to $x$ in the usual sense.
So $h_n:= \funk{\cdot}{x_n}{T}-\funk{b}{x_n}{T}$ converges pointwise
on $T$ to $h$ as $n\to\infty$.

By dropping initial terms if necessary, we may assume that $x_n$
is not greater than or equal to $y$ in the ordering on $T$ for any $n\in\N$.
So, for each $n\in\N$, $\funk{y}{\cdot}{T}$ is positive and increasing
along the straight line segment between $y$ and $x_n$.
Let $\gamma_n$ be this line segment parameterised in such a way that
$\gamma_n(0):=y$ and $\funk{y}{\gamma_n(t)}{T}:=t$
for all $0\le t\le \funk{y}{x_n}{T}$.

Fix $t>0$. Then, $\gamma_n(t)$ converges in the usual sense to
$\gamma(t)$ as $n\to\infty$.
By Lemma~\ref{lem:continuity}, this implies that each of
$\funk{\gamma(t)}{\gamma_n(t)}{T}$ and $\funk{\gamma_n(t)}{\gamma(t)}{T}$
converge to zero. Since $\Funky{T}$ satisfies the triangle inequality,
\begin{align*}
-\funk{\gamma(t)}{\gamma_n(t)}{T} \le h_n(\gamma_n(t)) - h_n(\gamma(t))
                                  \le \funk{\gamma_n(t)}{\gamma(t)}{T}.
\end{align*}
Using in addition the pointwise convergence of $h_n$, we conclude that
$h_n(\gamma_n(t))$ converges to $h(\gamma(t))$ as $n$ tends to infinity.

But $h_n(y)-h_n(\gamma_n(t)) = t$ for all $n\in\N$,
and so $h(y)-h(\gamma(t)) = t$.

Since this is true for all $t\in\rplus$, $\gamma$ is an almost--optimal path
with parameter $0$ with respect to the function $h$.
So by Lemma~3.4 of~\cite{walsh},
\begin{equation*}
h(y) \ge \h_{R,y}|_T(y)+\mumax_h(\h_{R,y}|_T).
\end{equation*}
Let $B:= \cB_T \intersection \closure\{g|_{T} \mid g\in \cK_{R} \}$.
From Lemma~3.6 of~\cite{AGW-m}, we know that $h\le w+\mumax_h(w)$
for all $w\in\closure \cK_T$.
So, since $\h_{R,y}|_T$ is in $B$,
which is a subset of $\closure \cK_T$ by Lemma~\ref{lem:lemmaB}, we have
\begin{equation*}
h(y) = \inf_{w\in B}
      \Big( w(y)+\mumax_h(w) \Big).
\end{equation*}
But $y$ is an arbitrary point of $T$, and so (\ref{eqn:representation}) holds
with $\nu:\cB_T\to\R\union\{+\infty\}$ defined by
\begin{equation*}
\nu(w) := \begin{cases}
   \mumax_h(w),
          & \text{if $w\in B$} \\
   +\infty, & \text{otherwise.}
\end{cases}
\end{equation*}
Since $B$ is closed in $\cB_T$
and $\mumax_h$ is lower semicontinuous, $\nu$ is lower semicontinuous.
\end{proof}

\begin{definition}
Let $T\subset V$ be an open cone, $x$ be in $\partial T$,
and $h$ be a function from $T$ to $\R$ satisfying
\begin{equation}
\label{homo}
h((1-\lambda)x + \lambda y) = \log \lambda + h(y)
\end{equation}
whenever $y$ and $(1-\lambda)x + \lambda y$ are in $T$ and $\lambda>0$.
We define the \emph{extension} of $h$ to $\tau(T,x)$ by
\begin{equation*}
h|^{\tau(T,x)}(y):= -\log \lambda_y + h((1-\lambda_y)x + \lambda_y y),
\qquad
\text{for all $y\in\tau(T,x),$}
\end{equation*}
where $\lambda_y>0$ is chosen so that $(1-\lambda_y)x + \lambda_y y$
is in $T$.
\end{definition}
Observe that the homogeneity condition~(\ref{homo}) implies that the
definition of $h|^{\tau(T,x)}(y)$ does not depend on the choice of $\lambda_y$.
By Lemma~\ref{lem:lemmaC}, for any open cone $T$ and point $p\in T$,
the function $f_{T,p}$ satisfies~(\ref{homo}) if $x\in[0]_T$.

\begin{lemma}
\label{lem:bigcone}
Let $T_1$ and $T_2$ be open cones in $V$ such that $T_1\subset T_2$.
Then $\funk{x}{y}{T_1} \ge \funk{x}{y}{T_2}$ for all $x,y\in T_1$.
\end{lemma}
\begin{proof}
If $\lambda y - x$ is in $T_1$ for some $\lambda>0$, then it is
also in $T_2$. Thus
\begin{equation*}
\{\lambda>0\mid \lambda y - x\in T_1\}
   \subset
\{\lambda>0\mid \lambda y - x\in T_2\}.
\end{equation*}
The result follows immediately.
\end{proof}

For the proof of the next lemma, we will need a tool from~\cite{AGW-m}.
Let $X$ be some set and $d(\cdot,\cdot)$ be a function on $X\times X$
satisfying the triangle inequality. We say a function $f$ from $X$ to $\R$
is $1$--Lipschitz if $f(x)-f(y)\le d(x,y)$ for all $x$ and $y$ in $X$.
This generalises the usual definition for metric spaces.
A function is said to be an \emph{extremal generator}
of the set of $1$--Lipschitz
functions if it is $1$--Lipschitz and can not be written as the minimum
of two $1$--Lipschitz functions each different from it.
It was shown in Theorem~6.2 of~\cite{AGW-m} that a function is an extremal
generator of the set of $1$--Lipschitz functions if and only if it is a
Busemann point or of the form $d(\cdot,x)-d(b,x)$ for some $x\in X$.

\begin{lemma}
\label{lem:extiter}
Let $T\subset V$ be an open cone and
let $\fun\in\cB_T\intersection A_T(x)$ with $x\in\partial T\backslash[0]_T$.
Then $\fun$ satisfies condition~(\ref{homo}) and so $\fun|^{\tau(T,x)}$,
its extension to $\tau(T,x)$, is well defined.
Moreover, $\fun|^{\tau(T,x)}\in \cK_{\tau(T,x)} \union \cB_{\tau(T,x)}$.
\end{lemma}
\begin{proof}
By Lemma~\ref{lem:representation}, $\fun$ can be expressed
as a min--plus combination of the elements of $\restsalt{T}{x}$, that is,
can be expressed in the form~(\ref{eqn:representation}).
Since each element of $\cK_{\tau(T,x)}$ satisfies the homogeneity
condition~(\ref{homo}), $\fun$ also satisfies it. Therefore the extension
$\fun|^{\tau(T,x)}$ is well defined.

Note that the extension of an infimum of functions satisfying~(\ref{homo})
is equal to the infimum of the extensions. It follows therefore from
Lemma~\ref{lem:representation} that $\fun|^{\tau(T,x)}$ can be written
as a min--plus combination of elements of $\closure \cK_{\tau(T,x)}$.
Since this set is the horofunction compactification of
$\Funky{\tau(T,x)}$, all its elements are 1-Lipschitz with respect
to $\Funky{\tau(T,x)}$.
It follows that $\fun|^{\tau(T,x)}$ is also 1-Lipschitz on $\tau(T,x)$
with respect to $\Funky{\tau(T,x)}$.

Now suppose that $\fun|^{\tau(T,x)}=\min(h_1,h_2)$,
where $h_1$ and $h_2$ are functions from $\tau(T,x)$ to $\R$ that are
1-Lipschitz in the Funk geometry on $\tau(T,x)$.
We deduce that $h_1|_T$ and $h_2|_T$ are 1-Lipschitz in the Funk
geometry on $T$ since, by Lemma~\ref{lem:bigcone},
$\Funky{T}$ is greater than $\Funky{\tau(T,x)}$ on $T$.

So, since $\fun$ is an extremal generator of the set of $1$--Lipschitz
functions of the Funk geometry on $T$, it must be equal to either
$h_1|_T$ or $h_2|_T$.
The homogeneity of $\fun|^{\tau(T,x)}$ and of $h_1$ and $h_2$ then
imply that $\fun|^{\tau(T,x)}$ equals either $h_1$ or $h_2$.
We have thus proved that $\fun|^{\tau(T,x)}$ is an extremal generator
of the set of $1$--Lipschitz functions of the Funk geometry on $\tau(T,x)$.
The conclusion follows.
\end{proof}

\begin{lemma}
\label{lem:busemannink}
Let $C\subset V$ be an open cone.
Every function in $\cK_C\union\cB_C$ can be written $\minf_{T,p}|_C$
for some $T\in\cones(C)$ and $p\in T$.
\end{lemma}
\begin{proof}
Let $h\in\cK_C\union\cB_C$.
If $h$ is in $\cK_C$, then $h=\minf_{C,p}$ for some $p$ in $C$.
Otherwise, $h\in\cB_C$ and since
$\bigcup_{x\in \partial C\backslash[0]_C} A_C(x) \supset \cB_C$,
we have that $h$ is in $A_C(x)$ for some $x$ in $\partial C\backslash[0]_C$.
By Lemma~\ref{lem:extiter},
$h|^{\tau(C,x)}$ is therefore in $\cK_{\tau(C,x)}\union\cB_{\tau(C,x)}$.
Using similar reasoning,
we deduce that $h|^{\tau(C,x)}$ is either in $\cK_{\tau(C,x)}$
or in $\cB_{\tau(C,x)} \intersection A_{\tau(C,x)}(y)$ for some
$y$ in $\partial \tau(C,x)\backslash[0]_{\tau(C,x)}$.
Proceeding inductively, we eventually reach some cone $T\in\cones(C)$
such that $h|^T$ is in $\cK_T$. The conclusion follows.
\end{proof}

\begin{lemma}
\label{lem:geoiter}
Let $T\subset V$ be an open cone and let $x\in \partial T\backslash[0]_T$.
Let $h\in \cK_{\tau(T,x)}\union \cB_{\tau(T,x)}$.
Then there exists a sequence of points in $T$ that is an almost-geodesic
with respect to
both the Funk and the reverse Funk metrics on $T$ and converges to $h|_T$
in the Funk sense and to $\hrev_{T,x}$ in the reverse Funk sense.
\end{lemma}
\begin{proof}
We write $R:=\tau(T,x)$ for convenience.
Let $(x_n)_{n\in\N}$ be an almost-geodesic in $R$ with respect to the
Funk metric on $R$ that converges to $h$ in the $\Funky{R}$ sense.
So $(x_n)_{n\in\N}$ satisfies, for some $\epsilon>0$,
\begin{equation}
\label{eqn:busemann1}
\sum_{n=0}^N \funk{x_n}{x_{n+1}}{R} \le \funk{x_0}{x_{N+1}}{R} + \epsilon,
\qquad
\text{for all $N\in\N$.}
\end{equation}

Let $S$ be a countable dense subset of $T$, and let $(z_n)_{n\in\N}$
be a sequence in $S$ that visits every point of $S$ infinitely often.
Let $y_n:=(1-\lambda_n)x + \lambda_n x_n$, where $(\lambda_n)_{n\in\N}$
is a sequence of positive real numbers which we have yet to specify.
We wish to choose $(\lambda_n)_{n\in\N}$ so that, for all $n\in\N$,
the following hold:
\begin{align}
\label{eqn:busemann2}
y_n &\in T \\
\label{eqn:busemann2bis}
|y_n-x| &< \frac{1}{n} \\
\label{eqn:busemann3}
\funk{z_n}{y_n}{T} - \funk{z_n}{y_n}{R} &< \frac{1}{n} \\
\label{eqn:busemann4}
\funk{b}{y_n}{T} - \funk{b}{y_n}{R} &< \frac{1}{n} \\
\label{eqn:busemann5}
\funk{y_{n}}{y_{n+1}}{T} - \funk{y_{n}}{y_{n+1}}{R} &< \frac{1}{2^n} \\
\label{eqn:busemann5bis}
\rev{y_{n}}{y_{n+1}}{T} + \hrev_{T,x}(y_{n+1}) - \hrev_{T,x}(y_{n})
    &< \frac{1}{2^n}.
\end{align}
Inclusion~(\ref{eqn:busemann2}) and inequality~(\ref{eqn:busemann2bis})
hold when $\lambda_n$ is small enough,
and, by Lemma~\ref{lem:lemmaB}, the same is true for~(\ref{eqn:busemann3})
and~(\ref{eqn:busemann4}).
On the other hand, inequalities~(\ref{eqn:busemann5})
and~(\ref{eqn:busemann5bis})
involve both $y_n$ and $y_{n+1}$, which means that the choice of
$\lambda_{n+1}$ must be made after that of $\lambda_n$.
So one must choose $\lambda_0,\lambda_1,\dots$ in that order.
That~(\ref{eqn:busemann5}) and~(\ref{eqn:busemann5bis}) may be satisfied once
$y_n$ has been fixed follows from, respectively, Lemma~\ref{lem:lemmaB}
and a combination of Lemmas~\ref{lem:revlemmaB} and~\ref{lem:revconv}.

By Lemma~\ref{lem:lemmaC},
\begin{align}
\label{eqn:busemann6}
\funk{y_n}{y_{n+1}}{R}
   &= \log\frac{\lambda_n}{\lambda_{n+1}} + \funk{x_n}{x_{n+1}}{R},
\qquad
\text{for all $n\in\N$, and} \\
\label{eqn:busemann7}
\funk{y_0}{y_{N+1}}{R}
   &= \log\frac{\lambda_0}{\lambda_{N+1}} + \funk{x_0}{x_{N+1}}{R},
\qquad
\text{for all $N\in\N$.}
\end{align}
Also, Lemma~\ref{lem:bigcone} gives that
\begin{equation}
\label{eqn:busemann8}
\funk{y_0}{y_{N+1}}{T} \ge \funk{y_0}{y_{N+1}}{R},
\qquad
\text{for all $N\in\N$.}
\end{equation}

To show that $(y_n)_{n\in\N}$ is an almost-geodesic with respect to the
Funk metric on $T$,
we vary $n$ in~(\ref{eqn:busemann5}) and~(\ref{eqn:busemann6})
between $0$ and $N$, and combine the resulting inequalities and
equalities with~(\ref{eqn:busemann1}),
(\ref{eqn:busemann7}), and~(\ref{eqn:busemann8}).
We get that, for all $N\in\N$,
\begin{equation*}
\sum_{n=0}^N \funk{y_n}{y_{n+1}}{T} - \funk{y_0}{y_{N+1}}{T}
   < \epsilon + \sum_{n=0}^N \frac{1}{2^n}
   < \epsilon + 2.
\end{equation*}
So $(y_n)_{n\in\N}$ is an almost-geodesic in the Funk metric on $T$
and must therefore converge in the $\Funky{T}$ sense.

Let $u\in S$.
By Lemma~\ref{lem:lemmaC},
\begin{align*}
\funk{u}{y_n}{R} - \funk{b}{y_n}{R} = \funk{u}{x_n}{R} - \funk{b}{x_n}{R},
\end{align*}
which by assumption converges to $h(u)$ as $n$ tends to infinity.
But we can find $n\in\N$ arbitrarily large such that $z_n=u$.
From inequalities~(\ref{eqn:busemann3}) and~(\ref{eqn:busemann4})
and Lemma~\ref{lem:bigcone},
\begin{equation*}
|\funk{u}{y_n}{T} - \funk{b}{y_n}{T} - \funk{u}{y_n}{R} + \funk{b}{y_n}{R} |
     < \frac{1}{n},
\end{equation*}
and, from the arbitrariness of $n$, we conclude that
$\funk{u}{y_n}{T} - \funk{b}{y_n}{T}$ converges to $h(u)$
as $n$ tends to $\infty$.
Since this holds for all $u$ in a dense subset of $T$,
we see that the almost-geodesic
$(y_n)_{n\in\N}$ converges to $h$ with respect
to the Funk metric on $T$.

We now wish to show that $(y_n)_{n\in\N}$ is an almost-geodesic
with respect to the reverse Funk metric on $T$.
We combine the inequalities obtained from~(\ref{eqn:busemann5bis})
by varying $n$ from  $0$ to $N$ to get
\begin{equation*}
\sum_{n=0}^N \rev{y_{n}}{y_{n+1}}{T} + \hrev_{T,x}(y_{N+1}) -\hrev_{T,x}(y_0)
   < \sum_{n=0}^N \frac{1}{2^n} < 2.
\end{equation*}
So $(y_n)_{n\in\N}$ is an almost--optimal path in the reverse Funk geometry
with respect to the 1--Lipschitz function $\hrev_{T,x}$.
It is therefore an almost-geodesic in the reverse Funk geometry by Lemma~7.1
of~\cite{AGW-m}.
By inequality~(\ref{eqn:busemann2bis}), $(y_n)_{n\in\N}$ converges in the
usual sense to $x$, and so by Lemma~\ref{lem:revconv}
it converges to $\hrev_{T,x}$ in the reverse Funk sense.
\end{proof}

\begin{lemma}
\label{lem:busemann}
Let $C\subset V$ be an open cone.
For each $T\in\cones(C)$ and $p\in T$,
the function $\minf_{T,p}|_C$ is in $\cK_C\union\cB_C$.
\end{lemma}
\begin{proof}
Since $T\in\cones(C)$, there exists a sequence of cones $(T_k)_{1\le k\le N}$
such that $T_k\in\Gamma(\{T_{k-1}\})$ for all $1<k\le N$,
and $T_1=C$ and $T_N=T$.
Observe also that $\minf_{T,p}\in\cK_T$.
It follows from Lemma~\ref{lem:geoiter} that if $S$ and $R$ are open cones
such that $R\in\Gamma(\{S\})$ and $h$ is a function in $\cK_R\union\cB_R$,
then $h|_S$ is in $\cK_S\union\cB_S$.
We apply this repeatedly to the situation above to deduce that
$\minf_{T,p}|_{T_n}$ is in $\cK_{T_n}\union\cB_{T_n}$ for all
$1\le n \le N$. Taking $n=1$ gives the result.
\end{proof}

\begin{prop}
\label{pro:funkbusemann}
Let $C$ be an open cone in $V$. Then
\begin{align*}
\cK_C\union\cB_C &= \Big\{ \minf_{T,p}|_C \mid \text{$T\in \cones(C)$
   and $p\in T$} \Big\}
\end{align*}
and the set of Busemann points of the Funk geometry on $C$ is precisely
\begin{align*}
\cB_C = \Big\{ \minf_{T,p}|_C \mid \text{$T\in \cones(C)\backslash \{C\}$
   and $p\in T$} \Big\}.
\end{align*}
\end{prop}
\begin{proof}
The first statement follows from Lemmas~\ref{lem:busemann}
and~\ref{lem:busemannink}, and the second follows immediately from this.
\end{proof}

Having determined the set of Busemann points of the Funk geometry,
we now investigate the question of when all horofunctions are Busemann points.
Our technique will be to reuse the key lemma in~\cite{walsh:normed} used
to answer the question for normed spaces.

Let $T\subset V$ be an open cone containing the basepoint $b$.
For each $x\in T$ and $y\in V$, define
\begin{align*}
\expf_{T,x}(y):= \frac{\M{y}{x}{T}}{\M{b}{x}{T}}.
\end{align*}
It follows immediately from Lemma~\ref{lem:dualformula} that $\expf_{T,x}$
is convex.
Note that $\expf_{T,x}(y)= \exp (\hfunk_{T,x}(y))$ for all $x$ and $y$ in $T$.

We will need some notions from convex analysis.
For a reference on this subject, the reader may consult~\cite{beer_book}.
We use the Painlev\'e--Kuratowski topology on the set of closed
sets of $V$.
In this topology, a sequence of closed sets $(C_n)_{n\in\N}$
is said to converge to a closed set $C$ if the upper and lower closed limits of
the sequence both equal $C$. These limits are defined to be, respectively,
\begin{align*}
\pul C_n &:= \bigcap_{n\ge 0} \closure \Big( \bigcup_{i>n} C_i \Big)
\qquad
\text{and} \\
\pll C_n &:=
   \bigcap\Big( \closure \bigcup_{i\ge 0} C_{n_i}
         \mid \text{$(n_i)_{i\in\N}$ is an increasing sequence in $\N$} \Big).
\end{align*}
An alternative characterisation of convergence is that $(C_n)_{n\in\N}$
converges to $C$ if and only if each of the following hold:
\begin{itemize}
\item
for each $x\in C$, there exists $x_n\in C_n$ for $n$ large enough, such that
$(x_n)_n$ converges to $x$.
\item
if $(C_{n_k})_{k\in\N}$ is a subsequence of the sequence of sets and
$x_k\in C_{n_k}$ for each $k\in\N$, then convergence of $(x_k)_{k\in\N}$
to $x$ implies that $x\in C$.
\end{itemize}

The Painlev\'e--Kuratowski topology can be used to define a topology on the
space of lower-semicontinuous functions as follows. Recall that the epigraph
of a function $f$ on $\R^n$ is the set
$\{(x,\alpha)\in \R^n\times\R\mid \alpha\ge f(x)\}$. A sequence of
lower-semicontinuous functions is declared to be convergent in the epigraph
topology if the associated epigraphs converge in the Painlev\'e--Kuratowski
topology on $\R^n\times\R$. For proper metric spaces,
the epigraph topology is identical to the Attouch--Wets topology.

\newcommand\indicator{I}

We use $\indicator_E$ to denote the indicator function, which takes
value $0$ on $E$ and $+\infty$ everywhere else.
The \emph{Legendre--Fenchel transform} of a function
$f:V \to \R\union\{\infty\}$ is the function $f^*:V^* \to \R\union\{+\infty\}$
defined by
\begin{equation*}
f^*(y):= \sup_{x\in V} \big( \dotprod{y}{x} - f(x) \big),
\qquad\text{for all $y\in V^*$}.
\end{equation*}
The Legendre--Fenchel transform is a bijection from the set of
proper lower-semicontinuous convex functions to itself
and is continuous in the epigraph topology.

Recall that the tangent cone of a cone $T$ at a point $x\in\partial T$
is defined to be
$\closure\{ \lambda(y-x)\mid \text{$\lambda>0$ and $y\in T$} \}$.
The following lemma relates this concept to the concept of
``open tangent cone'' defined earlier.
We use the notation $A^\circ$ to denote the polar of a set $A$.
\begin{lemma}
\label{lem:opendual}
Let $T$ and $C$ be open cones in $V$. Then $T$ is the open tangent cone to $C$
at some point $x$ if and only if $\closure T$ is the tangent cone to
$\closure C$ at $x$.
\end{lemma}
\begin{proof}
We first show that
\begin{align}
\label{eqn:opdual}
\closure\bigcup_{\lambda>0} {\lambda}(\closure C-x)
   = \closure\bigcup_{\lambda>0} {\lambda}(C-x).
\end{align}
That the right hand side is contained in the left follows from
$C\subset \closure C$. The opposite inclusion follows from the principle
that the union of closures is contained in the closure of the union.

It is immediately clear from~(\ref{eqn:opdual}) that if $T$ is the
open tangent cone to $C$ at $x$, then $\closure T$ is the tangent
cone to $\closure C$ at $x$.
The converse follows from~(\ref{eqn:opdual}) and the fact that both
$T$ and $\bigcup_{\lambda>0}{\lambda}(C-x)$ are open convex sets and
hence regular open, that is, each is equal to the interior of
its closure.
\end{proof}

Recall that an exposed face of a convex set is the intersection of the
set with a supporting hyperplane.
\begin{lemma}
\label{lem:dualexposed}
Let $C$ be a closed cone in $V$. A closed cone $T$ is a tangent cone to $C$
at some point in $\partial C$ if and only if the dual cone $T^*$
is an exposed face of $C^*$.
\end{lemma}
\begin{proof}
Let $T$ be the tangent cone to $C$ at $x\in\partial C$, that is
$T=\closure\bigcup_{\lambda>0} \lambda(C-x)$.
We have
\begin{align*}
T^* &= \bigcap_{\lambda>0} \frac{1}{\lambda}(C-x)^\circ \\
    &= \bigcap_{\lambda>0} \frac{1}{\lambda}(\conv(C \union \{-x\}))^\circ \\
    &= \bigcap_{\lambda>0} \frac{1}{\lambda}
                 (C^*\intersection\{z\in V^* \mid \dotprod{z}{x}\le 1 \} ) \\
    &= C^*\intersection\{z\in V^* \mid
               \text{$\dotprod{z}{x}\le 1/\lambda$ for all $\lambda>0$} \} \\
    &= C^*\intersection\{z\in V^* \mid \dotprod{z}{x}\le 0 \}.
\end{align*}
Recall that $x\in\partial C$ if and only if $\dotprod{z}{x}\ge 0$
for all $z\in C^*$ and $\dotprod{z}{x}= 0$ for some $z\in C^*$.
So $T^*$ is the intersection of $C^*$ with a supporting hyperplane and hence
is an exposed face of $C^*$.

The converse can be proved by reversing the argument.
\end{proof}

\begin{lemma}
\label{lem:dualextreme}
Let $C$ be an open cone in $V$. Then $T\in\cones(C)$ if and only if $T^*$
is an extreme set of $C^*$.
\end{lemma}
\begin{proof}
According to Lemma~3.2 of~\cite{walsh:normed}, $T^*$ is an extreme set
of $C^*$ if and only if there exists a finite sequence of convex sets
$F_0,\dots,F_n$ in $V^*$ such that $F_0=C^*$, $F_n=T^*$, and $F_{i+1}$
is an exposed face of $F_i$ for each $i\in\{0,\dots,n-1\}$.
The conclusion follows on combining this with Lemmas~\ref{lem:dualexposed}
and~\ref{lem:opendual}.
\end{proof}

\begin{lemma}
\label{lem:conjugates}
Let $T\subset V$ be an open cone.
For all $x\in T$, we have that $\expf^*_{T,x}$ is the indicator function of
the set
\begin{align}
\label{eqn:indicatorset}
Z_{T,x} := T^*\intersection\{z\in V^*\mid \M{b}{x}{T} \dotprod{z}{x} \le 1 \}.
\end{align}
\end{lemma}
\begin{proof}
Using the definition of $Z_{T,x}$, we calculate the Legendre--Fenchel
transform of $I_{Z_{T,x}}$:
\begin{align*}
I_{Z_{T,x}}^* (y)
   &= \sup\{\dotprod{w}{y}
             \mid \text{$w\in T^*$ and $\M{b}{x}{T}\dotprod{w}{x}\le 1$} \} \\
   &= \frac{1}{\M{b}{x}{T}} \sup\Big\{ \frac{\dotprod{w}{y}}{\dotprod{w}{x}}
          \mid w\in T^* \Big\}.
\end{align*}
But, by Lemma~\ref{lem:dualformula}, this last expression is equal to
$\expf_{T,x}(y)$. The conclusion follows on taking the
Legendre--Fenchel transform.
\end{proof}

\begin{lemma}
\label{lem:equilipschitzian}
Let $C$ be an open cone in $V$. Then the set of functions
$\{\expf_{C,x} \mid \text{$T\in\cones(C)$ and $x\in T$} \}$
is equi--Lipschitzian.
\end{lemma}
\begin{proof}
Let $T\in\cones(C)$. Then $T\supset C$, and so $T^*\subset C^*$.
By Lemma~\ref{lem:dualformula},
\begin{align*}
\topfn_C(b/x)\dotprod{z}{x} \ge \dotprod{z}{b},
\qquad
\text{for all $z\in C^*$.}
\end{align*}
We conclude that $Z_{C,x}$ is a subset of
$\{z\in C^* \mid \dotprod{z}{b}\le 1\}$.
Since we are assuming the base point $b$ is in the interior of $C$,
this set is compact.
So, by Lemma~\ref{lem:conjugates},
all the $\expf^*_{T,x}$; $T\in\cones(C)$, $x\in T$ take the value $+\infty$
outside the same bounded set.
It follows that their Legendre--Fenchel transforms are equi--Lipschitzian.
\end{proof}

\begin{lemma}
\label{lem:convergence}
Let $C\subset V$ be an open cone. Let $(T_n)_{n\in\N}$ be a sequence
of cones in $\cones(C)$ and let $(x_n)_{n\in\N}$ be a sequence
of points such that $x_n\in T_n$ for all $n\in\N$.
Then, $(\hfunk_{T_n,x_n}|_C)_{n\in\N}$ converges pointwise if and only if
$(\expf_{T_n,x_n})_{n\in\N}$ converges in the epigraph topology.
Moreover, the limit of the former sequence is $\hfunk_{T,x}|_C$,
for some $T\in\cones(C)$ and $x\in T$,
if and only if the limit of the latter is $\expf_{T,x}$.
\end{lemma}
\begin{proof}
Since, by Lemma~\ref{lem:equilipschitzian},
the $\expf_{T_n,x_n};\,n\in\N$ are equi--Lipschitzian, convergence of
$(\expf_{T_n,x_n})_{n\in\N}$ in the epigraph topology implies its
convergence pointwise~\cite[Prop.~7.1.3]{beer_book}, and therefore also
the pointwise convergence of $(\hfunk_{T_n,x_n}|_C)_{n\in\N}$.
Moreover, if $(\expf_{T_n,x_n})_{n\in\N}$ converges to $\expf_{T,x}$
with $T\in\cones(C)$ and $x\in T$, then the limit of
$(\hfunk_{T_n,x_n}|_C)_{n\in\N}$ will be $\hfunk_{T,x}|_C$.

Now suppose that $\hfunk_{T_n,x_n}|_C$ converges pointwise.
Then, $\expf_{T_n,x_n}$ converges pointwise on $C$. It follows from this
and the fact that the $\expf_{T_n,x_n}$ are equi--Lipschitzian that
$(g_{T_n,x_n})_{n\in\N}$ converges in the epigraph topology,
where
\begin{align*}
g_{T,x}: V\to\R,\, y\mapsto\begin{cases}
   \expf_{T,x}(y), & \text{if $y\in \closure C$}, \\
   +\infty,            & \text{otherwise},
      \end{cases}
\end{align*}
for all $T\in\cones(C)$ and $x\in T$.
So $g^*_{T_n,x_n}$ also converges in the same topology.

Since $g_{T,x} = \max\{ j_{T,x}, \indicator_{\closure C} \}$
for any cone $T\in\cones(C)$ and $x\in T$, we have
$g^*_{T,x} = \conv \{ \expf_{T,x}^* , I^*_{\closure C}\}$,
where $\conv$ denotes the convex hull of a set of functions.
By Lemma~\ref{lem:conjugates}, the Legendre--Fenchel transform of $j_{T,x}$
is $\indicator_{Z_{T,x}}$.
Also, since $C$ is a cone, the transform of $\indicator_{\closure C}$
is $\indicator_{-C^*}$. The convex hull of these two functions is
the indicator function of the convex hull of $Z_{T,x}$ and $-C^*$.
Since $-C^*$ is a cone, the convex hull of these two sets is equal to
\begin{align*}
Z_{T,x} - C^*
   := \left\{ y-w \in V^* \mid \text{$y\in Z_{T,x}$ and $w\in C^*$} \right\}.
\end{align*}

Let $z$ be in both $C^*$ and $Z_{T,x}-C^*$.
So $z=y-w$ for some $y\in Z_{T,x}$ and $w\in C^*$.
Therefore $y/2$ is in $T^*$ and is a convex combination of $z$ and $w$.
Since, by Lemma~\ref{lem:dualextreme}, $T^*$ is an extreme set of $C^*$,
we conclude that each of $z$ and $w$ are in $T^*$.
So $\dotprod{w}{x}\ge0$ and hence
\begin{align*}
\dotprod{z}{x} = \dotprod{y}{x} - \dotprod{w}{x}
               \le \dotprod{y}{x}
              \le \frac{1}{\M{b}{x}{T}}.
\end{align*}
We deduce that $z$ is in $Z_{T,x}$.
We have proved that
\begin{align*}
C^* \intersection (Z_{T,x} - C^*) \subset C^* \intersection Z_{T,x}.
\end{align*}
But the reverse inclusion is trivial since $0\in C^*$, and so we have the
equality of these two sets.

It follows that $g^*_{T,x}$ agrees with $\expf^*_{T,x}$ on $C^*$.

Since $\indicator_{Z_{T,x}}$ takes the value $+\infty$ outside $T^*$,
which is a subset
of $C^*$, we conclude that $\expf^*_{T,x}=\max(g^*_{T,x}, I_{C^*})$.

It follows from this and the convergence of $g^*_{T_n,x_n}$ in the epigraph
topology that $\expf^*_{T_n,x_n}$, and hence $\expf_{T_n,x_n}$, converges
in the same topology.

Running through the same argument assuming that $\hfunk_{T_n,x_n}|_C$
converges to $\hfunk_{T,x}|_C$ for some $T\in\cones(C)$ and $x\in T$,
we get that $\expf_{T_n,x_n}$ converges to $\expf_{T,x}$.
\end{proof}

Define
\begin{align*}
\Lambda: \R^n\times (0,1] \to \R^n\times \R_+,\,
   (x,y)\mapsto \Big(\frac{x}{y}, \frac{1}{y}-1\Big).
\end{align*}

\begin{lemma}
\label{lem:collineation}
The map $\Lambda$ is a homeomorphism. Furthermore, each of $\Lambda$ and
$\Lambda^{-1}$ map line segments to line segments.
\end{lemma}
\begin{proof}
The inverse of $\Lambda$ is
\begin{align*}
\Lambda^{-1}: \R^n\times \R_+ \to \R^n\times (0,1],\,
   (x,y)\mapsto \Big(\frac{x}{1+y}, \frac{1}{1+y}\Big).
\end{align*}
So $\Lambda$ is a bijection and clearly each of $\Lambda$ and $\Lambda^{-1}$
is continuous.

The second part of the lemma follows from the fact that $\Lambda$ is a
linear--fractional function~\cite{boyd_vandenberghe}.
\end{proof}

For any open cone $C\subset V$, define
\begin{align*}
\zsets_C := \Big\{Z_{T,x}\backslash\{0\}
                  \mid \text{$T\in\cones(C)$ and $x\in T$}\Big\}.
\end{align*}

\begin{lemma}
\label{lem:Zclosed}
Let $C\subset V$ be an open cone.
Then, $\cK_C\union\cB_C$ is closed in the pointwise topology if and only
if $\zsets_C$ is closed in the Painlev\'e--Kuratowski topology.
\end{lemma}
\begin{proof}
The map assigning to each closed set
of a metric space its indicator function is an embedding of the non-empty
closed sets with the Painlev\'e--Kuratowski topology into the proper
lower--semicontinuous functions with the epigraph
topology~\cite[Proposition 7.1.1]{beer_book}.
Combining this with the fact that the Legendre--Fenchel transform
is continuous in the epigraph topology, we see from Lemma~\ref{lem:conjugates}
that $\zsets_C$ is closed in the Painlev\'e--Kuratowski topology
if and only if the set of functions
$J_C := \{ \expf_{T,x} \mid \text{$T\in\cones(C)$ and $x\in T$} \}$
is closed in the epigraph topology.

Assume $\cK_C\union \cB_C$ is closed and let $(\expf_{T_n,x_n})_{n\in\N}$
be a sequence of functions in $J_C$ convergent in the epigraph topology.
Then $f_{T_n,x_n}|_C$ converges pointwise by Lemma~\ref{lem:convergence}.
Since $f_{T_n,x_n}|_C$ is in $\cK_C\union \cB_C$ for all $n\in\N$
by Proposition~\ref{pro:funkbusemann},
the limit must be in $\cK_C\union \cB_C$.
So, by Lemma~\ref{lem:busemannink}, we can write the limit as $f_{T,x}|_C$
for some $T\in\cones(C)$ and $x\in T$. Using Lemma~\ref{lem:convergence} again,
we see that $(\expf_{T_n,x_n})_{n\in\N}$ converges to $\expf_{T,x}$.
Therefore $J_C$ is closed.

The converse may be proved in a similar manner.
\end{proof}

Choose a coordinate system on $V$ so that $b=(0,\dots,0,1)$.
Using these coordinates, we may consider the sets $\R^{N-1}\times(0,1]$
and $\R^{N-1}\times\rplus$ to be subsets of $V^*$, where $N$ is the dimension
of $V$.
The map $\Lambda$ is a bijection between these sets, and maps
$\{z\in C^*\backslash\{0\}\mid\dotprod{z}{b}\le1\}$ to $B^\circ\times\rplus$,
where $B^\circ$ is the subset of $\R^{N-1}$ such that
$B^\circ\times\{1\}=\{z\in C^*\mid\dotprod{z}{b}=1\}$.
We can extend $\Lambda$ in an obvious way to a bijection $\tilde\Lambda$
between subsets
of $\R^{N-1}\times(0,1]$ and subsets of $\R^{N-1}\times\rplus$.

Let $\dualdbps$ be the set of functions from $\R^{N-1}$ to
$\R\union\{+\infty\}$ that are finite and affine on an extreme set of
$B^\circ$, take the value $+\infty$ everywhere else, and have infimum zero.
Let $\dbps:=\{ f^* \mid f\in \dualdbps \}$ be the set of Legendre--Fenchel
transforms of these functions and let $U_C$ be the set of their epigraphs.

\begin{lemma}
\label{lem:Uclosed}
Let $C\subset V$ be an open cone. Then,
$\zsets_C$ is closed in the Painlev\'e--Kuratowski topology
if and only if $U_C$ is.
\end{lemma}
\begin{proof}
Let $T\in \cones(C)$ and $x\in T$.
By Lemma~\ref{lem:dualformula}, $\M{b}{x}{T}\dotprod{z}{x}\ge\dotprod{z}{b}$
for all $z\in T^*$, and so $Z_{T,x}$ is a subset of
$\{z\in C^*\mid\dotprod{z}{b}\le1\}$.
We see also that $Z_{T,x}$ is the intersection of the cone $T^*$
and a half-space. So, by Lemma~\ref{lem:collineation},
$\Omega := \tilde\Lambda(Z_{T,x}\backslash\{0\})$ is the intersection
of the vertical cylinder $E\times\R$ and a half-space,
where $E$ is the subset of $\R^{N-1}$ such that
$E\times\{1\}=\{z\in T^*\mid\dotprod{z}{b}=1\}$.
So we may think of $\Omega$ as the epigraph
of a function $f$ that is affine on $E$ and takes the value
$+\infty$ outside this set.
By Lemma~\ref{lem:dualextreme}, $E$ is an extreme set of $B^\circ$.
We have that
\begin{align*}
\inf f &= \inf\{t\in\rplus \mid
               \text{$(p,t)\in\Omega$ with $p\in \R^{N-1}$} \} \\
       &= \inf\left\{\frac{1}{\dotprod{z}{b}}-1
                      \mid z\in Z_{T,x}\backslash\{0\}\right\} \\
       &= \frac{1}{\sup\{\dotprod{z}{b}
           \mid \text{$z\in T^*$ and $\M{b}{x}{T}\dotprod{z}{x}\le 1$}\}} -1 \\
       &= \Big(\sup_{z\in T^*}
           \frac{\dotprod{z}{b}}{\M{b}{x}{T}\dotprod{z}{x}}\Big)^{-1} -1.
\end{align*}
The last line equals zero by Lemma~\ref{lem:dualformula}.

Therefore $f$ is in $\dualdbps$ and so $\Omega$ is in $U_C$.
We conclude that $\tilde\Lambda$ maps sets in $\zsets_C$ to sets in $U_C$.

That $\tilde\Lambda^{-1}$ maps sets in $U_C$ to sets in $\zsets_C$ may be
established in a similar manner.

Observe that each element of $\zsets_C$ is closed in $\R^n\times (0,1]$
and each element of $U_C$ is closed in $\R^n\times \R_+$.
Since $\Lambda$ is a homeomorphism, $\tilde\Lambda$ is a homeomorphism
between the set of closed
subsets of $\R^n\times (0,1]$ and those of $\R^n\times \R_+$, the
Painlev\'e--Kuratowski topology being used in both cases.
Therefore $\zsets_C$ is closed in this topology if and only if $U_C$ is.
\end{proof}

\begin{prop}
\label{pro:funk3}
Let $C$ be an open cone in $V$.
A necessary and sufficient condition for every horofunction of the Funk
geometry on $C$ to be a
Busemann point is that the set of extreme sets of the dual cone $C^*$
be closed in the Painlev\'e--Kuratowski topology.
\end{prop}
\begin{proof}
First assume that every horofunction in the Funk geometry is a Busemann point,
in other words that $\cK_C\union\cB_C=\closure\cK_C$.
Then $\cK_C\union\cB_C$ is closed and so, by Lemma~\ref{lem:Zclosed},
$\zsets_C$ is closed in the Painlev\'e--Kuratowski topology.
We now use Lemma~\ref{lem:Uclosed} to deduce that $U_C$ is also closed
in this topology.
Since the elements of $U_C$ are the epigraphs of the elements of $\dualdbps$,
we conclude that $\dualdbps$ is closed in the epigraph topology.
It follows that $\dbps$, the set of Legendre--Fenchel transforms of
these functions is also closed in this topology.
But in~\cite[Lemma~4.3]{walsh:normed}, it was shown that this is equivalent
to the set of extreme sets of $B^\circ$ being closed in the
Painlev\'e--Kuratowski topology, and this is equivalent to the closure
of the set of extreme sets of $C^*$.

To establish the reverse implication, we assume that the set of extreme sets
of $C^*$ is closed and reverse the chain of argument.
We conclude that $\cK_C\union\cB_C$ is closed.
Since $\cK_C \subset \cK_C\union\cB_C \subset \closure\cK_C$,
it follows that $\cK_C\union\cB_C = \closure\cK_C$,
which implies that every horofunction is a Busemann point.
\end{proof}

\section{Boundary of the Hilbert geometry}

Determining the horofunction boundary of the Hilbert geometry will involve
combining what we know about the boundaries of the Funk and reverse Funk
geometries.

\begin{lemma}
\label{lem:splitseq}
A sequence in an open cone $C$ is an almost-geodesic in the Hilbert geometry
if and only if it is an almost-geodesic in both the Funk and reverse Funk
geometries.
\end{lemma}
\begin{proof}
Let $(x_n)_{n\in\N}$ be a sequence in $C$. For all $n\in\N$, define
\begin{align*}
L(n) &:= \sum_{i=1}^n \hil{x_{i-1}}{x_i}{C} - \hil{x_0}{x_n}{C} \\
F(n) &:= \sum_{i=1}^n \funk{x_{i-1}}{x_i}{C} - \funk{x_0}{x_n}{C} \\
R(n) &:= \sum_{i=1}^n \rev{x_{i-1}}{x_i}{C} - \rev{x_0}{x_n}{C}.
\end{align*}
Then, for all $n\in\N$, we have that $L(n)=F(n)+R(n)$
and that $L(n)$, $F(n)$, and $R(n)$
are all non-negative.

If $(x_n)_{n\in\N}$ is an almost-geodesic in the Hilbert geometry,
then, for some $\epsilon>0$,
we have that $L(n)<\epsilon$ for all $n\in\N$.
It follows that $R(n)$ and $F(n)$ are also less than $\epsilon$
for all $n\in\N$, which means that $(x_n)_{n\in\N}$ is an almost-geodesic
in  the Funk and reverse Funk geometries.

On the other hand, if, for some $\epsilon_1>0$ and $\epsilon_2>0$,
we have $F(n)<\epsilon_1$ and $R(n)<\epsilon_2$ for all $n\in\N$,
then $L(n)<\epsilon_1+\epsilon_2$ for all $n\in\N$,
proving that $(x_n)_{n\in\N}$ is an almost-geodesic in the
Hilbert geometry.
\end{proof}

Recall that $A_C(x)$ is the set of Funk-metric
horofunctions that may be attained as a limit in the Funk sense of a
sequence converging to $x\in \partial C\backslash[0]_C$ in the usual sense.

\begin{lemma}
\label{lem:hilberthorofunctions}
The set of horofunctions of the Hilbert geometry on an open cone $C$
containing no lines is
\begin{equation*}
\{  \hrev_{C,x} + \hfunk \mid 
   \text{$x\in\euclidean$ and $\hfunk\in A_C(x)$} \}.
\end{equation*}
\end{lemma}
\begin{proof}
If $\hfunk\in A_C(x)$, then there exists a sequence $(x_n)_{n\in\N}$ in $C$
converging in the usual sense to $x$ and in the Funk sense to $\hfunk$.
So
\begin{equation*}
\hil{\cdot}{x_n}{C}-\hil{b}{x_n}{C}
   = \rev{\cdot}{x_n}{C}-\rev{b}{x_n}{C}
         + \funk{\cdot}{x_n}{C}-\funk{b}{x_n}{C} 
\end{equation*}
converges to $\hrev_{C,x} + \hfunk$ by Lemma~\ref{lem:revconv}.
This proves that $\hrev_{C,x} + \hfunk$ is a horofunction of the
Hilbert geometry.

Now suppose that a sequence $(x_n)_{n\in\N}$ in $C$ converges in the Hilbert
sense to a horofunction $\hhil$. We may assume that $(x_n)_{n\in\N}$ is
contained in a cross section $D$ of the cone and,
since $C$ contains no lines, we may take $\closure D$ to be compact.
So some subsequence $(x_{n_i})_{i\in\N}$ converges in the usual sense
to some point $x$ in $\closure D$.
But $x$ cannot be in $D$ since $\hhil$ would then be equal to
$\hil{\cdot}{x}{C}-\hil{b}{x}{C}$ and we have assumed that $\hhil$
is a horofunction. Therefore $x$ is in the relative boundary of $D$.
So $\rev{\cdot}{x_{n_i}}{C}-\rev{b}{x_{n_i}}{C}$ converges to
$\hrev_{C,x}$ by Lemma~\ref{lem:revconv},
and therefore $\funk{\cdot}{x_{n_i}}{C}-\funk{b}{x_{n_i}}{C}$
converges to a function $\hfunk\in A_C(x)$ such that
$\hhil = \hrev_{C,x} + \hfunk$.
\end{proof}

Recall that $\cB_C$ is the set of Busemann points of the Funk geometry
on $C$.
For each open cone $C$ and $x\in\partial C\backslash\{0\}$,
define the set of functions
\begin{equation*}
B(x):=\{ \hrev_{C,x}+\hfunk\mid 
   \text{$\hfunk\in \cB_C \intersection A_C(x)$} \}.
\end{equation*}

\begin{lemma}
\label{lem:extremalhilbert}
Let $C$ be an open cone containing no lines.
The set of Busemann points of the Hilbert geometry on $C$ is
$\minspacehilbert := \bigcup_{x\in\euclidean} B(x)$.
If $h\in B(x)$ with $x\in \euclidean$,
then there exists an almost-geodesic converging to $h$ in
the Hilbert sense that converges to $x$ in the usual sense.
\end{lemma}
\begin{proof}
Assume that $h$ is a Busemann point.
So there exists an almost-geodesic
$(x_n)_{n\in\N}$ converging to $h$ in the Hilbert sense.
We may assume that this almost-geodesic lies in a cross section $D$
of the cone.
By Lemma~\ref{lem:splitseq}, the sequence $(x_n)_{n\in\N}$ is also
an almost-geodesic in both the Funk and reverse Funk geometries.
So, by Proposition~7.3 of~\cite{AGW-m}, it converges to a Busemann
point in both geometries; denote the respective limits by
$\hfunk$ and $\hrev$.
We see from Proposition~\ref{pro:reversehorofunctions} that
$r=\hrev_{C,x}$ for some $x$ in the relative boundary of $D$, and that
$(x_n)_{n\in\N}$ must converge in the usual sense to $x$.
Therefore, $\hfunk\in \cB_C \cap A_C(x)$, and so $h=\hrev_{C,x}+f$ is in $B(x)$,
and therefore in $\minspacehilbert$.

Conversely, assume that $\hfunk\in \cB_C \intersection A_C(x)$.
Then, by Lemma~\ref{lem:extiter}, $\hfunk|^{\tau(C,x)}$,
its extension to $\tau(C,x)$, is in $\cK_{\tau(C,x)} \union \cB_{\tau(C,x)}$.
Applying Lemma~\ref{lem:geoiter},
we obtain an almost-geodesic in the Funk geometry on $C$ converging
in the Funk sense to $\hfunk$ that is also an almost-geodesic in
the reverse Funk geometry, with respect to which it converges to
$\hrev_{C,x}$.
We may assume that this almost-geodesic lies in a cross section $D$
of the cone.
By Lemma~\ref{lem:splitseq}, this path must also be an almost-geodesic
in the Hilbert geometry on $C$, converging to $\hrev_{C,x}+\hfunk$.
Therefore, this function is a Busemann point of the
Hilbert geometry.
This concludes the proof of the first statement.

To prove the second, it suffices to observe that the convergence of the
almost--geodesic in the reverse Funk sense to $\hrev_{C,x}$ that we
have just established implies by Proposition~\ref{pro:reversehorofunctions}
its convergence to a positive multiple of $x$. By rescaling, we can
make the almost--geodesic converge to $x$ itself.
\end{proof}

\begin{proof}[Proof of Theorem~\ref{theorem1}]
Let $x\in \partial C \backslash\{0\}$.
Lemmas~\ref{lem:extiter} and~\ref{lem:geoiter}
and Proposition~\ref{pro:reversehorofunctions} together imply that
\begin{align*}
\cB_{C} \intersection A_C(x)
   = \Big\{ g|_C \mid g\in\cK_{\tau(C,x)} \union \cB_{\tau(C,x)} \Big\}.
\end{align*}
But, by Proposition~\ref{pro:funkbusemann},
\begin{align*}
\cK_{\tau(C,x)} \union \cB_{\tau(C,x)}
   = \Big\{\minf_{T,p}|_{\tau(C,x)} \mid
          \text{ $T\in\cones(\tau(C,x))$ and $p\in T$} \Big\}.
\end{align*}
The conclusion now follows from Lemma~\ref{lem:extremalhilbert}.
\end{proof}

\begin{lemma}
\label{lem:disjoint}
Let $C$ be an open cone containing no lines.
If $x$ and $y$ are on distinct rays of $\partial C$, then
$B(x)$ and $B(y)$ are disjoint.
\end{lemma}
\begin{proof}
Assume that $h$ is in both $B(x)$ and $B(y)$, with $x$ and $y$
lying on distinct rays of $\partial C$.
By Lemma~\ref{lem:extremalhilbert},
there exist almost-geodesics $(x_n)_{n\in\N}$ and
$(y_n)_{n\in\N}$ each converging to $h$ in the Hilbert sense
that converge respectively to $x$ and $y$ in the usual sense.
We may assume that $x$ and $y$ and both almost--geodesics
lie in some cross section $D$ of the cone $C$.

We define another sequence of points inductively in the following way.
Choose an initial point $z_0\in C$ arbitrarily. For $n$ odd, choose
$z_n:= x_m$ for some $m\ge n$ large enough that
\begin{equation}
\label{eqn:disjoint1}
\hil{z_{n-1}}{z_n}{C} + h(z_n) \le h(z_{n-1}) + \frac{1}{2^n}.
\end{equation}
This will be possible since, by Lemma~3.3 of~\cite{walsh},
\begin{equation*}
\lim_{m\to\infty} \Big(\hil{z_{n-1}}{x_m}{C} + h(x_m) \Big) = h(z_{n-1}).
\end{equation*}
Similarly, for $n$ even, we choose $z_n:= y_m$ for $m\ge n$
large enough that~(\ref{eqn:disjoint1}) holds.

We vary $n$ in~(\ref{eqn:disjoint1}) from $1$ to any $N\ge 1$ and
combine the resulting inequalities. We obtain that, for all $N\ge 1$,
\begin{equation*}
-h(z_0) + \sum_{n=1}^N \hil{z_{n-1}}{z_n}{C} + h(z_{N})
   \le \sum_{n=1}^N \frac{1}{2^n}
   < 1.
\end{equation*}
In other words, $(z_n)_{n\in\N}$ is an almost-optimal path with respect to $h$
with parameter~$1$ in the Hilbert geometry.
Therefore, by Lemma~7.1 of~\cite{AGW-m}, it is an almost-geodesic
in this geometry.
By Lemma~\ref{lem:splitseq}, it must also be an almost-geodesic in
the reverse Funk geometry, and so,
by Proposition~\ref{pro:reversehorofunctions}, it must converge in the
usual sense to some point in the relative boundary of $D$.
However, it contains infinitely many terms of each of the sequences
$(x_n)_{n\in\N}$ and $(y_n)_{n\in\N}$, and therefore both $x$ and $y$ are
limit points of $(z_n)_{n\in\N}$ in the usual topology. It follows that $x=y$.
\end{proof}

\newcommand\gfunk{q}
\begin{lemma}
\label{thm:uniquesplit}
Let $C\subset V$ be an open cone containing no lines.
Every Hilbert--geometry horofunction $h$ can be written in a unique
way as $h = \hrev_{C,x} + \gfunk$, with $x\in \partial C\backslash\{0\}$
and $\gfunk\in A_C(x)$.
\end{lemma}
\begin{proof}
That $h$ can be written in the above form was proved in
Lemma~\ref{lem:hilberthorofunctions}.

Suppose that
\begin{equation*}
h= \hrev_{C,x} + \gfunk_{x} = \hrev_{C,y} + \gfunk_{y}
\end{equation*}
with $\gfunk_{x}\in A_C(x)$ and $\gfunk_{y}\in A_C(y)$.
By Lemma~\ref{lem:representation}, we may find a min-plus measure $\mu$ on
$\cB_C$ such that
\begin{equation*}
h= \inf_{w\in\cB_C} \Big(\hrev_{C,x} + w + \mu(w) \Big)
\end{equation*}
and $\mu$ takes the value $+\infty$ outside $\restsalt{C}{x}$.
By Lemma~\ref{lem:lemmaB}, $\restsalt{C}{x} \subset A_C(x)$.
Lemma~\ref{lem:extremalhilbert} says that
$\hrev_{C,x} + w$ is a Busemann point of the Hilbert geometry
for all $w\in \cB_C\intersection A_C(x)$.
We deduce that
\begin{equation*}
h= \inf_{z\in\minspacehilbert} \Big( z + \bar\mu(z) \Big),
\end{equation*}
where $\bar\mu:\minspacehilbert\to \R\union\{+\infty\}$ is defined by
\begin{equation*}
\bar\mu(z):=
   \begin{cases}
   \mu(z-\hrev_{C,x}), & \text{if $z\in B(x)$}, \\
   +\infty, & \text{otherwise}.
   \end{cases}
\end{equation*}
Observe that $\bar\mu$ is lower semicontinuous since
$\mu$ is lower semicontinuous on $\cB_C$ and $B(x)$ is closed
in $\minspacehilbert$.
Therefore $\bar\mu$ is a min-plus measure on $\minspacehilbert$.

Using similar reasoning, we may find a min-plus measure $\bar\nu$
on $\minspacehilbert$ taking the value $+\infty$ outside $B(y)$
such that $h=\inf_{z\in\minspacehilbert} (z + \bar\nu(z))$.

By Theorem~1.1 of~\cite{walsh}, there exists a min-plus measure $\rho$ on
$\minspacehilbert$ satisfying $h=\inf_{z\in\minspacehilbert} (z + \rho(z))$
that is greater than or equal to each of
$\bar\mu$ and $\bar\nu$.
Since $h$ is not identically $+\infty$, neither is $\rho$,
and therefore $B(x)$ and $B(y)$ must have an element in common.
So, by Lemma~\ref{lem:disjoint}, $x$ is a positive multiple of $y$.
Therefore, $\hrev_{C,x}=\hrev_{C,y}$,
and hence $\gfunk_{x} = \gfunk_{y}$.
\end{proof}

\begin{proof}[Proof of Theorem~\ref{theorem3}]
Let $C$ be an open cone in $\R^{N+1}$ of which $D$ is a cross section.
Since $D$ is bounded, $C$ contains no lines.

Comparing Lemmas~\ref{lem:hilberthorofunctions} and~\ref{lem:extremalhilbert},
we see that if all horofunctions of the Funk geometry on $C$
are Busemann points then so also are all horofunctions of the
Hilbert geometry on $C$.

Now assume that all Hilbert horofunctions are Busemann points
and let $f$ be a Funk geometry horofunction.
So $f\in A_C(x)$ for some $x\in \partial C\backslash\{0\}$.
By Lemma~\ref{lem:hilberthorofunctions}, $\hrev_{C,x}+ f$ is a Hilbert geometry
horofunction, and so, by Lemma~\ref{lem:extremalhilbert},
can be written as $\hrev_{C,y}+ g$ for some $y\in\partial C\backslash\{0\}$
and $g\in \cB_C\intersection A_C(y)$.
But from Lemma~\ref{thm:uniquesplit}, we must have that $y$ is a multiple
of $x$ and that $f=g$. Therefore, $f$ is a Busemann point of the
Funk geometry.

We have shown that all Hilbert horofunctions are Busemann points if and
only if all Funk horofunctions are. Applying Proposition~\ref{pro:funk3},
we see that this in turn is equivalent to the closure of the set of
extreme sets of $C^*$ in the Painlev\'e--Kuratowski topology.
But the polar of $D$ may be identified with a cross section of $C^*$
and so its set of extreme sets is closed in the Painlev\'e--Kuratowski
topology if and only if the set of extreme sets of $C^*$ is.
\end{proof}

\begin{proof}[Proof of Theorem~\ref{theorem2}]
Let $(x_n)_{n\in\N}$ be a sequence in $D$ that converges in the Hilbert
sense to a point $h$ in the horofunction boundary of the Hilbert geometry.
Let $y$ and $z$ in $\partial D$ be limit points of $(x_n)_{n\in\N}$
in the usual topology,
and let $(y_n)_{n\in\N}$ and $(z_n)_{n\in\N}$ be subsequences of
$(x_n)_{n\in\N}$ converging in the usual sense to $y$ and $z$ respectively.
By Lemma~\ref{lem:revconv},
$(y_n)_{n\in \N}$ converges to $\hrev_{C,y}$
and $(z_n)_{n\in \N}$ converges to $\hrev_{C,z}$ in the reverse Funk sense,
where $C$ is some open cone in $\R^{N+1}$ of which $D$ is
a cross section.
It follows that $y_n$ converges in the Funk sense to
$\gfunk_{y}:=h-\hrev_{C,y}$, which must be in $A_C(y)$,
and similarly $z_n$ converges in the Funk sense
to $\gfunk_{z}:=h-\hrev_{C,z}$, which must be in $A_C(z)$. So
\begin{equation*}
h=\hrev_{C,y} + \gfunk_{y} = \hrev_{C,z} + \gfunk_{z}.
\end{equation*}
We now apply Lemma~\ref{thm:uniquesplit} to conclude that
\begin{equation*}
\hrev_{C,y} = \hrev_{C,z}
\qquad\text{and}\qquad
\gfunk_{y} = \gfunk_{z}.
\end{equation*}
The first of these implies that $y=z$ by
Proposition~\ref{pro:reversehorofunctions}.
So the sequence $(x_n)_{n\in\N}$ has only one limit point in $\closure D$.
This implies that it has a limit since $\closure D$ is compact.
\end{proof}

\section{Examples}
\label{sec:examples}
                                                                                
\begin{figure}
\label{fig:horospheres}
\centering
\begin{tabular}{cc}
{\includegraphics{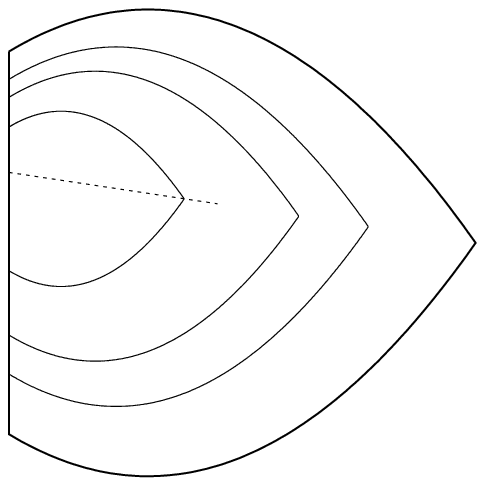}}
                       & \hspace{1.5cm}{\includegraphics{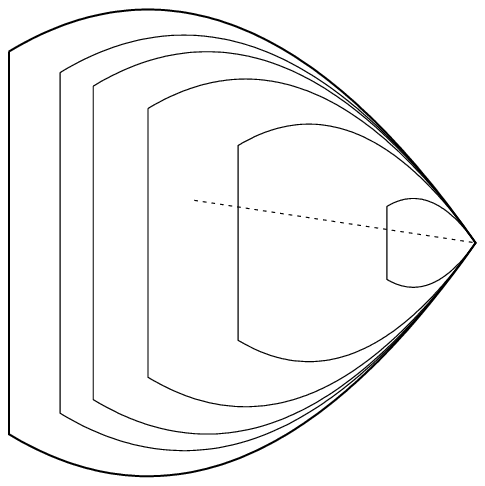}} \\
\includegraphics{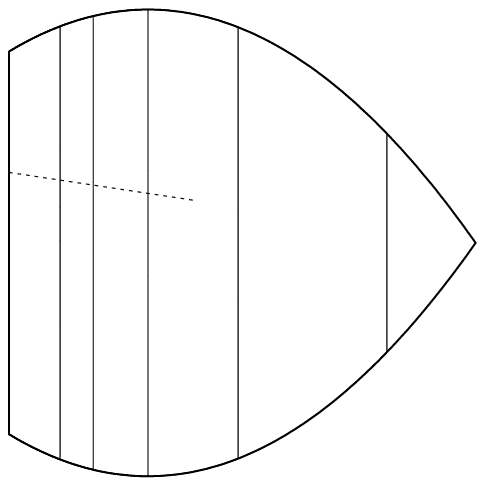}
                       & \hspace{1.5cm}\includegraphics{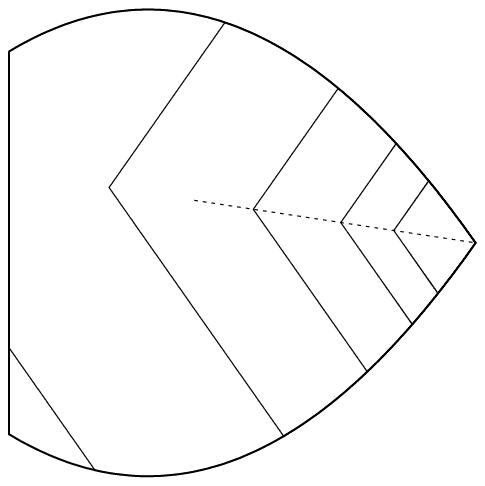} \\
{\includegraphics{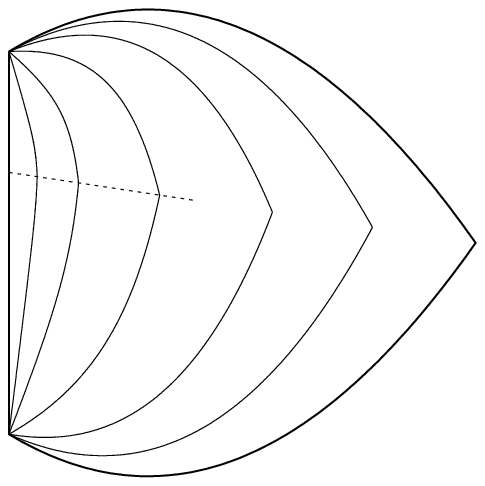}}
                       & \hspace{1.5cm}{\includegraphics{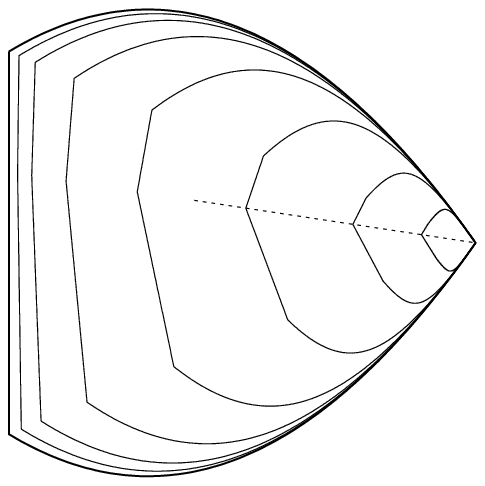}}
\end{tabular}
\caption{Horospheres in the reverse-Funk (top), Funk (middle) and
Hilbert (bottom) geometries, associated to two different geodesics.}
\end{figure}

In Figure~\ref{fig:horospheres}, we see some examples of horofunctions
in the reverse-Funk, Funk and Hilbert geometries associated to a
particular 2-dimensional convex domain.
Plotted are the level sets of the horofunctions (the horospheres).
On the left is, for each geometry, the limiting horofunction of a
straight-line geodesic (dotted) approaching the boundary at a flat point.
On the right are similar plots when the geodesic approaches a
non-differentiable point of the boundary.

In connection with Theorem~\ref{theorem3}, we give some examples of domains
where the set of extreme sets is closed and some where it is not.

\begin{example}
In dimension two, the set of extreme sets of any convex set is always
closed. Therefore, horofunctions of a Hilbert geometry
are always Busemann points.
\end{example}
                                                                                
\begin{example}
In dimension three, define the set
\begin{align*}
D := \Big\{ (x,y,z)\in\R^3
             \mid \text{$|x|+|z|\le 1$ and $x^2+y^2\le 1$} \Big\}.
\end{align*}
The polar of $D$ is the convex hull of the square with corners
$(\pm1,0,\pm1)$ and the circle $\{(x,y,z)\mid \text{$x^2+y^2=1$, $z=0$}\}$.
                                                                                
For all $n\in\N$, let $p_n:=(\cos (1/n),\,\sin (1/n),\,0)$.
Observe that the sequence of extreme sets $(\{p_n\})_{n\in\N}$
of $D^\circ$ converges to the set $\{(1,0,0)\}$ as $n\to\infty$.
However, this set is not extreme.
                                                                                
So from Theorem~\ref{theorem3} we would expect the existence of a
horofunction that is not a Busemann point.
One can show that the function $f:D\to\R,\,(x,y,z)\mapsto \log(1-x)$
is a non-Busemann horofunction of the Funk geometry on $D$.
Adding the reverse-Funk horofunction associated to the point $(1,0,0)$,
one obtains a  non-Busemann horofunction of the Hilbert geometry.
\end{example}
                                                                                
\begin{example}
In dimension three, the set of extreme sets of a convex set is closed
if and only if the set of extreme points is. So, for an example showing that
closure of the set of extreme points of the dual ball is not sufficient for all
horofunctions to be Busemann points, one must go to dimension four.
                                                                                
Let $D$ be the polar of the closed convex hull of the four circles
\begin{align*}
S_1^{\pm} &:= \Big\{ (x,y,\pm 1, 0)\in\R^4
             \mid \text{$x^2 + y^2 = 1$} \Big\} \\
S_2^{\pm} &:= \Big\{ (\pm 1, 0, w, z)\in\R^4
             \mid \text{$w^2 + z^2 = 1$} \Big\}.
\end{align*}
It was shown in~\cite{walsh:normed} that the set of extreme points of $D^\circ$
is closed but the set of its extreme sets is not.
So one would expect here also non-Busemann horofunctions of the Funk and
Hilbert geometries.
\end{example}

\bibliographystyle{plain}
\bibliography{hilbert}

\end{document}